\documentclass[11pt]{amsart}
\usepackage{amsmath,amssymb,epsfig,color, amsthm, mathrsfs}
\usepackage{wrapfig,lipsum,floatrow,sidecap,comment}
\usepackage{graphicx}
\usepackage{multirow}
\usepackage{hhline}
\usepackage{url}
\usepackage{mathtools}
\graphicspath{ {images/} }
\allowdisplaybreaks

\graphicspath{ {images/} }

\def\b{\textbf{b}}
\def\L{\mathcal{L}}

\def\p{\prime}

\newtheorem{theorem}{Theorem}[section]

\newtheorem{definition}[theorem]{Definition}

\newtheorem{proposition}[theorem]{Proposition}
\newtheorem{remark}[theorem]{Remark}
\newtheorem{example}[theorem]{Example}

\begin{document}
\title{The Braid Indices of Pretzel Links: \\
A Comprehensive Study, Part I}
\author{Yuanan Diao$^\dagger$, Claus Ernst$^*$ and Gabor Hetyei$^\dagger$}
\address{$^\dagger$ Department of Mathematics and Statistics\\
University of North Carolina Charlotte\\
Charlotte, NC 28223}
\address{$^*$ Department of Mathematics\\
Western Kentucky University\\
Bowling Green, KY 42101, USA}
\email{}
\subjclass[2010]{Primary: 5725; Secondary: 5727}
\keywords{knots, links, braid index, alternating links, Seifert graph, Morton-Frank-Williams inequality, pretzel link, Montesinos link.}

\begin{abstract}
The determination of the braid index of an oriented link is generally a hard problem. In the case of alternating links, some significant progresses have been made in recent years which made explicit and precise braid index computations possible for links from various families of alternating links, including the family of all alternating Montesinos links. However, much less is known for non-alternating links. For example, even for the non-alternating pretzel links, which are special (and simpler) Montesinos links, the braid index is only known for a very limited few special cases. In this paper and its sequel, we study the braid indices for all non-alternating pretzel links by a systematic approach. We classify the pretzel links into three different types according to the Seifert circle decompositions of their standard link diagrams. More specifically, if $D$ is a standard diagram of an oriented pretzel link $\L$, $S(D)$ is the Seifert circle decomposition of $D$, and $C_1$, $C_2$ are the Seifert circles in $S(D)$ containing the top and bottom long strands of $D$ respectively, then $\L$ is classified as a Type 1 (Type 2) pretzel link if $C_1\not=C_2$  and $C_1$, $C_2$ have different (identical) orientations. In the case that $C_1=C_2$, then $\L$ is classified as a Type 3 pretzel link. In this paper, we present the results of our study on Type 1 and Type 2 pretzel links. Our results allow us to  determine the precise braid index for any non-alternating Type 1 or Type 2 pretzel link. Since the braid indices are already known for all alternating pretzel links from our previous work, it means that we have now completely determined the braid indices for all Type 1 and Type 2 pretzel links.
\end{abstract}

\maketitle
\section{Introduction}\label{sec:intro}

A well known and important invariant in knot theory is the braid index of a link, which is defined as the minimum number of strings in a braid needed to represent the said link as a closed braid. The determination of the braid index of a link used to be an extremely difficult problem. In fact, prior to the discovery of the HOMFLY-PT polynomial~\cite{HOMFLY,PT}, the braid index was known for very few knots and links. The HOMFLY-PT polynomial, more specifically the Morton-Frank-Williams inequality (MFW inequality) derived from it~\cite{FW,Mo}, provided us a powerful tool for braid index determination. In the case of  
alternating links, precise braid index formulas have been derived for various link families including the two bridge links~\cite{Mu}, the alternating pretzel links~\cite{Diao2021,MP} and more generally the family of all alternating Montesinos links~\cite{Diao2021}. In the case of non-alternating links, the braid index formulas are known for all closed positive braids with a full twist (which include all the torus links)~\cite{FW} and some special link families such as those defined in~\cite{Lee} by the so called pretzel graphs. However, in general, we know much less about the braid indices of non-alternating links. For example, despite their relatively simple structures, the braid indices of non-alternating pretzel links remain unsettled, with very few exceptions. This fact indicates that there are additional challenges when dealing with non-alternating links, and it is likely some new techniques/methods in addition to those used for the alternating links are needed. 

\medskip
Motivated by the problem outlined above, in this paper and its sequel, we study the braid index problem for all  pretzel links in a systematic and thorough manner. Since the braid index problem has already been solved for the alternating pretzel links in \cite{Diao2021}, the new challenge comes from the non-alternating pretzel links. We are happy to report that we are able to determine the precise braid index for the vast majority of pretzel links. The main approach here is similar to that used in~\cite{Diao2021, MP} and is based on two results. The first one, which is due  to Yamada~\cite{Ya}, states that the number of Seifert circles in a link diagram $D$, denoted by $s(D)$ throughout this paper, is an upper bound for the braid index $\b(D)$ of $D$. The second result is the Morton-Frank-Williams inequality (MFW inequality)~\cite{FW,Mo}, which states that the braid index $\b(\L)$ of any link $\L$ is bounded below by $(E(\L)-e(\L))/2+1$, which we shall denote by $\b_0(\L)$, where $E(\L)-e(\L)$ is the degree span of the variable $a$ in the HOMFLY-PT polynomial $H(\L,z,a)$ of $\L$ as defined in Definition \ref{Def1} (with $E(\L)$ and $e(\L)$ being the highest and lowest powers of the variable $a$ in $H(\L,z,a)$). This approach succeeds if we can find a diagram $D$ of $\L$ such that $s(D)=\b_0(\L)$, since it then leads to $s(D)=\b_0(\L)\le \textbf{b}(\L)\le s(D)$ hence $\b(\L)=\b_0(\L)$. Although the approach sounds simple and straight forward, there is no guarantee that it would succeed since we may or may not be able to determine $E(\L)$ and $e(\L)$ for the given link $\L$. And even if we succeed with the first step, we still may not be able to construct a link diagram $D$ of $\L$ with $s(D)=\b_0(\L)=(E(\L)-e(\L))/2+1$. When this happens, it is difficult to know whether it is because there are no such diagrams exist, or it is because we have not found the right way to construct such a diagram as we know that the MFW inequality is strict for some links (for example there are five knots with crossing number less than or equal to 10 for which $\b(\L)>\b_0(\L)$~\cite{Adams}). These obstacles are the reasons why only partial results are known for the non-alternating pretzel links up to date. The braid index formulas (that we obtain in this paper) for some non-alternating pretzel links resemble those of their alternating counterparts, but for some other non-alternating pretzel links, are very different from those of their alternating counterparts, where we have to develop new methods/techniques to overcome these new challenges. In the following, let us first introduce the concept of pretzel links and how we classify them into different types.

\begin{figure}[htb!]
\includegraphics[scale=.8]{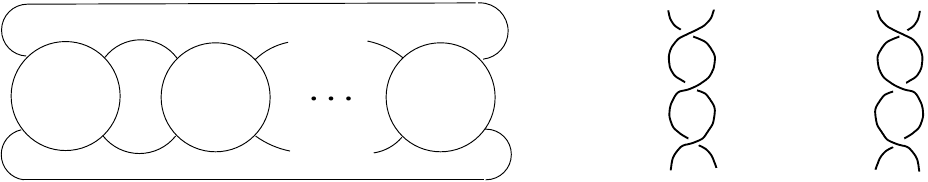}
\caption{Left: The general structure of a Montesinos link: each circle in the picture contains a (2 string) rational tangle; Middle: a vertical strip of 2 right-handed half twists; Right: a vertical strip of 2 left-handed half twists.
}
\label{Mont} 
\end{figure}

\medskip
A Montesinos link is a link that has a diagram with a structure as shown in Figure \ref{Mont}, in which each circle contains a rational tangle.  A pretzel link is a special Montesinos link where each rational tangle consists of only a vertical strip of (at least one) right handed or left handed half twists. For each pretzel link we consider in this paper, it is understood that an orientation has been assigned to each of its components.  

\medskip

\medskip
\begin{definition}\label{Def2}{\em 
Let $D$ be a pretzel link diagram and $S(D)$ be the Seifert circle decomposition of $D$. Let $C_1$ and $C_2$ be the Seifert circles in $S(D)$ containing the top and bottom strands of $D$ respectively. Then $D$ is said to be a Type 1 (2) pretzel link if $C_1\not=C_2$  and $C_1$, $C_2$ have different (same) orientations. In the case that $C_1=C_2$, then $D$ is said to be a Type 3 pretzel link.}
\end{definition}

\medskip
\begin{figure}[htb!]
\includegraphics[scale=.7]{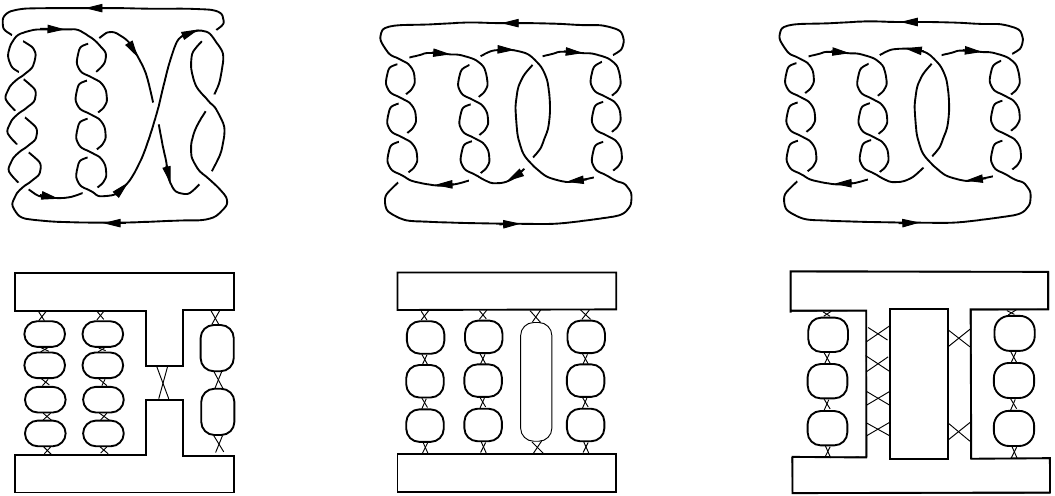}
\caption{From left to right: Examples of Type 1, 2 and 3 pretzel link diagrams and their corresponding Seifert circle decompositions. 
}
\label{Types} 
\end{figure}

Figure \ref{Types} shows an example of pretzel link for each of Types 1, 2 and 3. Notice that $D$ is of Type 1 (2) if and only every strip in $D$ contains an odd (even) number of crossings (half twists) and these crossings all smooth horizontally. We note that these types of pretzel links are classified as Types M1, M2 and B Montesinos links in~\cite{Diao2021}. Also, as the examples in Figure \ref{Types} show, different orientation assignments to the same underline diagram can lead to different types of pretzel links. 
Without loss of generality, we can always orient the top long strand in a pretzel link diagram from right to left. This choice, together with the crossing sign and the information on how the crossings in a strip are to be smoothed (either all vertically or all horizontally), will allow us to determine the handedness of all crossings in the diagram. Thus in our notations introduced below, we only need to indicate the sign of the crossings in each strip  (where the sign of a crossing is defined in Figure \ref{fig:cross}) and how the crossings are to be smoothed. 

\medskip
\begin{definition}\label{TypeDef}{\em 
We introduce the following {\em grouping} of the pretzel links as follows.
We denote by $P_{1}(2\alpha_1+1,\ldots, 2\alpha_{\kappa^+}+1; -(2\beta_1+1),\ldots, -(2\beta_{\kappa^-}+1))$ the set of all Type 1 pretzel links with $\kappa^+$ strips containing $2\alpha_1+1,\ldots, 2\alpha_{\kappa^+}+1$ positive crossings respectively and $\kappa^-$ strips containing $2\beta_1+1,\ldots, 2\beta_{\kappa^-}+1$ negative crossings respectively. Similarly, $P_{2}(2\alpha_1,\ldots, 2\alpha_{\kappa^+}; -2\beta_1,\ldots, -2\beta_{\kappa^-})$ denotes  the set of all Type 2 pretzel links with $\kappa^+$ strips containing $2\alpha_1,\ldots, 2\alpha_{\kappa^+}$ positive crossings respectively and $\kappa^-$ strips containing $2\beta_1,\ldots, 2\beta_{\kappa^-}$ negative crossings respectively. For a Type 3 pretzel link diagram, we first divide its strips of half-twists into two parts. If the crossings in a strip are smoothed vertically, then the strip is placed into Part 1, otherwise the strip is placed into Part 2. $P_3(\mu_1,\ldots,\mu_{\rho^+};-\nu_1,\ldots,-\nu_{\rho^-}\vert 2\alpha_1,\ldots, 2\alpha_{\kappa^+}; -2\beta_1,\ldots, -2\beta_{\kappa^-})$ denotes the set of all Type 3 pretzel link diagrams in which $\mu_1,\ldots,\mu_{\rho^+}$ ($\nu_1,\ldots,\nu_{\rho^-}$) are the numbers of crossings in the positive (negative) strips in Part 1, and $2\alpha_1,\ldots, 2\alpha_{\kappa^+}$ ($2\beta_1,\ldots, 2\beta_{\kappa^-}$) are the numbers of crossings in the positive (negative) strips in Part 2. Notice that in the case of a Type 3 pretzel link, we must have $\kappa^++\kappa^-=2n$ for some $n\ge 1$, due to the fact that two adjacent Seifert circles sharing crossings must have opposite orientations unless they are concentric to each other. For the sake of simplicity, we do not allow one crossing strips with different crossing signs in a pretzel link diagram since such crossings can be pairwise deleted via Reidemeister II moves. That is, we shall assume that in our definition here that either $\alpha_j>0$ for all $j$, or $\beta_i>0$ for all $i$ in a Type 1 pretzel link, and either $\mu_j>1$ for all $j$ or $\nu_i>1$ for all $i$ in a Type 3 pretzel link. 
}
\end{definition}

\medskip
A pretzel link diagram with its top strand oriented from right to left and satisfying the condition of Definition \ref{TypeDef} is called a {\em standard} pretzel link diagram in this paper. Since the union of the sets of all Types 1, 2 and 3 as defined in Definition \ref{TypeDef} is the entire set of all pretzel links (although they do not form a partition), every pretzel link can be represented by a standard diagram. For the sake of convenience, throughout this paper, we will only be using standard pretzel link diagrams. We should point out that two pretzel links from the same $P_1$, $P_2$ or $P_3$ set as defined in Definition \ref{TypeDef} may or may not be topologically equivalent. However, as our results will show, any two pretzel links from the same set have the same braid index. We summarize our results in the following theorems. The fact that there are many different cases to consider in our results is a testimony on the complexity of this problem.

\medskip
\begin{theorem}\label{MT1}
If $\L\in P_{1}(2\alpha_1+1,\ldots, 2\alpha_{\kappa^+}+1; -(2\beta_1+1),\ldots, -(2\beta_{\kappa^-}+1))$, then the braid index of $\L$, denoted by $\b(\L)$, is given by the following formulas:
\begin{equation}\label{MT1e1}
\b(\L)=\left\{
\begin{array}{ll}
\sum \alpha_j\ &{\rm if}\ \kappa^+=2, \kappa^-=1, \beta_1=0,\\
\sum \beta_i\ &{\rm if}\ \kappa^+=1, \kappa^-=2, \alpha_1=0,
\end{array}
\right.
\end{equation}
\begin{equation}\label{MT1e2}
\b(\L)=\left\{
\begin{array}{ll}
1+\sum \alpha_j\ &{\rm if}\ 
\left\{
\begin{array}{l}
\kappa^+=2, \kappa^-=1, 0<\beta_1<\min\{\alpha_1,\alpha_2\},\ {\rm or} \\
\kappa^+\ge 3, \kappa^-=1, \beta_1<\min_{1\le j\le \kappa^+}\{\alpha_j\},
\end{array}
\right.\\
1+\sum \beta_i\ &{\rm if}\ 
\left\{
\begin{array}{l}
\kappa^-=2, \kappa^+=1, 0<\alpha_1<\min\{\beta_1,\beta_2\},\ {\rm or} \\
\kappa^-\ge 3, \kappa^+=1, \alpha_1<\min_{1\le i\le \kappa^-}\{\beta_i\},
\end{array}
\right.\\
\end{array}
\right.
\end{equation}
\begin{equation}\label{MT1e3}
\b(\L)=\left\{
\begin{array}{ll}
2+\sum \alpha_j+ \beta_1-\min_{1\le j\le \kappa^+}\{\alpha_j\}\ &{\rm if}\ 
\ \kappa^+\ge 2, \kappa^-=1, \beta_1\ge \min_{1\le j\le \kappa^+}\{\alpha_j\}\\
2+\sum \beta_i+ \alpha_1-\min_{1\le i\le \kappa^-}\{\beta_i\}\ &{\rm if}\ 
\ \kappa^-\ge 2, \kappa^+=1, \alpha_1\ge \min_{1\le i\le \kappa^-}\{\beta_i\},
\end{array}
\right.
\end{equation}
\begin{equation}\label{MT1e4}
\b(\L)=1+\sum \alpha_j+\sum \beta_i\ {\rm if}\ \tau=\kappa^+-1\ {\rm or}\ \tau^\prime=\kappa^--1
\end{equation}
and
\begin{equation}\label{MT1e5}
\b(\L)=2+\sum \alpha_j+\sum \beta_i\ {\rm \quad for \ all\ other\ cases},
\end{equation}
where $\tau$ ($\tau^\prime$) is the number of $\beta_i$'s ($\alpha_j$'s) that equal to zero and it is understood that $\sum \alpha_j=0$ if $\kappa^+=0$ and $\sum \beta_i=0$ if $\kappa^-=0$ in {\em (\ref{MT1e5})}.
\end{theorem}

We note that the braid index formula for the Type 1 alternating pretzel links, namely the special cases $\kappa^+=0$ or $\kappa^-=0$ in (\ref{MT1e5}), has been obtained in \cite{Diao2021, MP}.  All other results in Theorem \ref{MT1} are new to our knowledge.

\medskip
\begin{theorem}\label{MT2}
If $\L\in P_{2}(2\alpha_1,\ldots, 2\alpha_{\kappa^+}; -2\beta_1,\ldots, -2\beta_{\kappa^-})$, then 
\begin{equation}\label{MT2e1}
\b(\L)=\left\{
\begin{array}{ll}
1+\sum \alpha_j&{\rm if}\ \kappa^-=0, \\
1+\sum \beta_i&{\rm if}\  \kappa^+=0,
\end{array}
\right.
\end{equation}
\begin{equation}\label{MT2e2}
\b(\L)=\left\{
\begin{array}{ll}
\sum \alpha_j&{\rm if}\ \kappa^+\ge 2, \kappa^-=1, \beta_1<\min_{1\le j\le \kappa^+}\{\alpha_j\},\\
\sum \beta_i&{\rm if}\ \kappa^-\ge 2, \kappa^+=1, \alpha_1<\min_{1\le i\le \kappa^-}\{\beta_i\},\\
\end{array}
\right.
\end{equation}
\begin{equation}\label{MT2e3}
\b(\L)=\left\{
\begin{array}{ll}
1+\sum \alpha_j+\beta_1-\min_{1\le j\le \kappa^+}\{\alpha_j\},&{\rm if}\ \kappa^+\ge 2, \kappa^-=1, \beta_1\ge \min_{1\le j\le \kappa^+}\{\alpha_j\},\\
1+\sum \beta_i+\alpha_1-\min_{1\le i\le \kappa^-}\{\beta_i\},&{\rm if}\ \kappa^-\ge 2, \kappa^+=1, \alpha_1\ge \min_{1\le i\le \kappa^-}\{\beta_i\},\\
\end{array}
\right.
\end{equation}
and 
\begin{equation}\label{MT2e4}
\b(\L)=\sum \alpha_j+\sum \beta_i\ {\rm if}\ \kappa^-\ge 2, \kappa^+\ge 2.
\end{equation}
\end{theorem}

\medskip
We note that (\ref{MT2e1}), namely the braid index formula for the Type 2 alternating pretzel links, has been obtained in \cite{Diao2021, MP}. We  include it here for the sake of completeness. Furthermore, some partial results for other cases have also been established in \cite{MP}. More specifically, in the cases of (\ref{MT2e2}) and (\ref{MT2e3}), it has been shown in \cite{MP} that $\b(\L)\le 1+\sum \alpha_j+\sum \beta_i$, and
in the case of  $\kappa^-\ge 2$ and $\kappa^+\ge 2$, it has been shown in \cite{MP} that $\sum \alpha_j+\sum \beta_i\le \b(\L)\le 1+\sum \alpha_j+\sum \beta_i$.

\medskip
\begin{example}{\em  The pretzel links in $P_1(5,3;-5,-1)$ (one diagram from this set is shown at the left in Figure \ref{Types}) have braid index  $1+\sum \alpha_j+\sum \beta_i=1+2+1+2+0=6$ by (\ref{MT1e4}). In this particular case, the alternating counterpart of $P_1(5,3;-5,-1)$ is  $P_1(5,5,3,1;0)$, whose pretzel links have braid index 7 by \cite[Theorem 4.7]{Diao2021}. The pretzel links in $P_2(4,4,2;-4)$ (one diagram from this set is shown in the center of Figure \ref{Types}) have 
braid index $1+\sum \alpha_j+\beta_1-\min\{\alpha_j\}=7$ by (\ref{MT2e3}). The alternating counterpart of $P_2(4,4,2;-4)$ is $P_2(4,4,4,2;0)$, whose pretzel links have braid index 8 by \cite[Theorem 4.7]{Diao2021}. }
\end{example}

\smallskip
\begin{remark}{\em We note that the difference between the braid indices of a Type 1 non alternating pretzel link $\L$ and its alternating counterpart can be one of the following cases: zero if $\L$ falls into the cases of (\ref{MT1e5}), 1 if $\L$ falls into the cases of (\ref{MT1e4}), 2 if $\L$ falls into the cases of (\ref{MT1e1}) and this difference can be any positive integer if $\L$ falls into the cases of (\ref{MT1e2}) or (\ref{MT1e3}). On the other hand, the difference between the braid indices of a Type 2 non alternating pretzel link $\L$ and its alternating counterpart is either 1, if $\L$ falls into the cases of (\ref{MT2e4}), or can be any positive integer if $\L$ falls into the cases of (\ref{MT2e2}) or (\ref{MT2e3}). (\ref{MT2e1}) does not apply since the pretzel links in this formula are alternating.
}
\end{remark}

\medskip
We will prove Theorems \ref{MT1} and \ref{MT2} in two parts. In the first part we establish that the formulas in the theorems are braid index lower bounds for the corresponding pretzel links. We will do this using the MFW inequality in Sections  \ref{lowerbound1} and \ref{lowerbound2}. This, of course, depends heavily on the HOMFLY-PT polynomials of these pretzel links and we will need to introduce some known results concerning the HOMFLY-PT polynomials. We shall do this in Section \ref{sec:prep}. In the second part of the proof, we will establish that the formulas in the theorems are braid index upper bounds for the corresponding pretzel links. This is done by direct constructions. We shall demonstrate these constructions in Sections \ref{upperbound1} and  \ref{upperbound2}. Some well known construction methods, as well as some new ones developed in this paper, will be introduced prior to that in Section \ref{upper_prep}. 

\section{preparations for establishing the lower bounds}\label{sec:prep}

\begin{figure}[htb!]
\includegraphics[scale=.6]{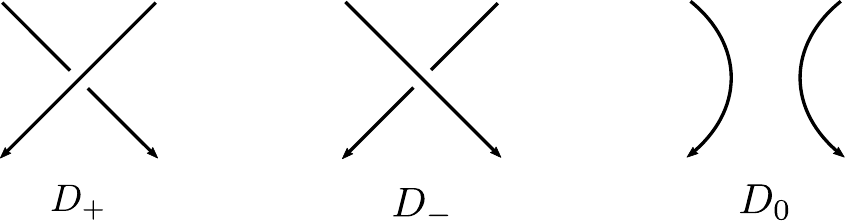}
\caption{The sign convention at a crossing of an oriented link and the splitting of the crossing: the crossing in $D_+$ ($D_-$) is positive (negative) and is assigned $+1$ ($-1$) in the calculation of the writhe of the link diagram.}
\label{fig:cross}
\end{figure}

\begin{definition}\label{Def1}{\em 
 Let $D$ be a link diagram of an oriented link $L$. Let $D_+$, $D_-$ and $D_0$ be the three link diagrams identical to $D$ except at one crossing as shown in Figure \ref{fig:cross}, then we define the HOMFLY-PT polynomial $H(D,z,a)$ by the following conditions:
\begin{enumerate}
\item If $D_1$ and $D_2$ are ambient isotopic, then $H(D_1,z,a)=H(D_2,z,a)$;
\item $aH(D_+,z,a) - a^{-1}H(D_-,z,a) = zH(D_0,z,a)$;
\item If $D$ is an unknot, then $H(D,z,a)=1$. 
\end{enumerate}
By~\cite{HOMFLY,Ja,PT}, these conditions uniquely define a polynomial $H(L,z,a)$ of $L$ independent of $D$.}
\end{definition}

\begin{figure}[htb!]
\includegraphics[scale=.75]{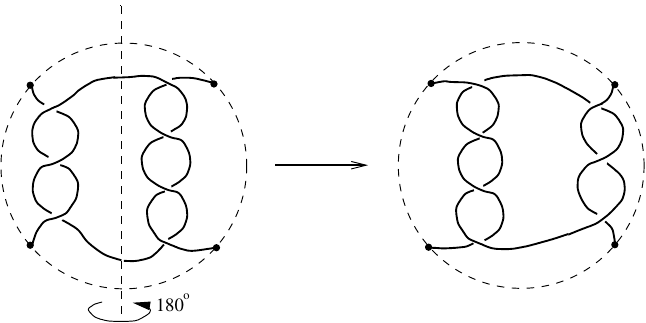}
\caption{A mutation move applied to two adjacent vertical strips in a pretzel link diagram.}
\label{mutation_move}
\end{figure}

\medskip
\begin{remark}\label{default} {\em 
A mutation move is an operation on a link diagram as shown in Figure \ref{mutation_move}. It is known that the HOMFLY-PT polynomial does no change under a mutation move~\cite[Proposition 11]{LICKORISH}. One can verify that a mutation move involving two adjacent vertical strips of an oriented pretzel link diagram will not change the crossing signs of the strips nor will it change how the crossings in the strip are smoothed (vertically or horizontally). It follows that changing the order of the crossing strips in a pretzel link will not change the HOMFLY-PT polynomial of the link. Therefore,  if
$\L\in P_1(2\alpha_1+1,\ldots,2\alpha_{\kappa^+}+1;-(2\beta_1+1),\ldots,-(2\beta_{\kappa^-}+1))$, $\L^\p\in P_1(2\alpha_1^\p+1,\ldots,2\alpha^\p_{\kappa^+}+1;-(2\beta^\p_1+1),\ldots,-(2\beta^\p_{\kappa^-}+1))$, then $H(\L)=H(\L^\p)$ as long as $\{\alpha_1,\ldots,\alpha_{\kappa^+}\}=\{\alpha^\p_1,\ldots,\alpha^\p_{\kappa^+}\}$ and $\{\beta_1,\ldots,\beta_{\kappa^+}\}=\{\beta^\p_1,\ldots,\beta^\p_{\kappa^+}\}$. The same statement holds for Type 2 and Type 3 pretzel links.
}
\end{remark}

\medskip
For our purpose, it is more convenient to use the following two equivalent forms of the skein relation as defined in Definition \ref{Def1}.
\begin{eqnarray}
H(D_+,z,a)&=&a^{-2}H(D_-,z,a)+a^{-1}zH(D_0,z,a),\label{Skein1}\\
H(D_-,z,a)&=&a^2 H(D_+,z,a)-azH(D_0,z,a).\label{Skein2}
\end{eqnarray}
For a link diagram $D$, we shall use $c(D)$ ($c^-(D)$) to denote the number of crossings (negative crossings) in $D$, and $w(D)=c(D)-2c^-(D)$ to denote the writhe of $D$. If we write $H(D,z,a)$ as a Laurent polynomial of $a$, then we will use  $p^h(D)$ ($p^\ell(D)$) to denote the coefficient of $a^{E(D)}$ ($a^{e(D)}$). $p^h(D)$ and $p^\ell(D)$ are Laurent polynomials of $z$ and we shall use $p^h_0(D)$ and $p^\ell_0(D)$ to denote the highest power terms of $p^h(D)$ and $p^\ell(D)$ respectively. This means that $p^h_0(D)$ and $p^\ell_0(D)$ are monomials of $z$.

\medskip
\begin{remark}\label{VP}{\em 
For a link diagram $D$, if we cannot determine $E(D)$ and $e(D)$ by direct computation, we will have to employ either (\ref{Skein1}) or (\ref{Skein2}). We start with a suitably chosen crossing. Say the crossing is positive and the skein relation (\ref{Skein1}) is applied at this crossing. Assuming that we are able to determine $E(D_-)$ and $E(D_0)$, then we compare  $-2+E(D_-)$ and $-1+E(D_0)$. If $-2+E(D_-)<-1+E(D_0)$, then $E(D)=-1+E(D_0)$ and $p^h(D)=zp^h(D_0)$. If $-2+E(D_-)>-1+E(D_0)$, then $E(D)=-2+E(D_-)$ and $p^h(D)=p^h(D_-)$. If $-2+E(D_-)=-1+E(D_0)$, then we have to compare $p^h(D_-)$ and $zp^h(D_0)$. As long as we can establish that $p^h(D_-)+zp^h(D_0)\not=0$, we can still conclude that $E(D)=-1+E(D_0)=-2+E(D_-)$. For example we would know $p^h(D_-)+zp^h(D_0)\not=0$ if $p^h_0(D_0)$ and $p^h_0(D_-)$ have the same sign.
Similarly, we can use this approach to determine $e(D)$ and possibly $p^\ell(D)$ or $p^\ell_0(D)$. This procedure will be used repeatedly in the next section since our proofs are almost always induction based. For the sake of convenience, we shall follow~\cite{Diao2021} and call this procedure the {\em Verification Procedure plus (minus)} (VP$^+$ or VP$^-$ for short) when a positive (negative) crossing is used for this procedure. Furthermore, when there are apparent Reidemeister I or II moves in $D_-$, $D_+$ or $D_0$ when VP$^+$ or VP$^-$ is applied, we shall assume such moves will be automatically done in our subsequent calculations concerning these diagrams, even though we will still be using the same notation for the sake of simplicity.
}
\end{remark}

\medskip
The following results concerning the HOMFLY-PT polynomial of a connected sum of torus links with one or two components will be needed in Sections \ref{lowerbound1} and \ref{lowerbound2}. These can be derived from the HOMFLY-PT polynomials of torus links (which are well known and easily computed directly) and the fact that $H(\L_1\#\L_2)=H(\L_1)H(\L_2)$ in general.

\medskip
\begin{remark} \label{toruslink_formula} 
{\em
Let $T_o(2k, 2)$ and $T_o(-2k,2)$ be the torus links whose components have opposite orientations and whose crossings are positive and negative respectively, then   
\begin{eqnarray}
H(T_o(2k,2))&=&\left\{
\begin{array}{l}
z(a^{-1} +a^{-3}+\ldots +a^{-2k+3})+(z+z^{-1})a^{-2k+1} -z^{-1}a^{-2k-1}, \ k>1,\\
(z+z^{-1})a^{-1} -z^{-1}a^{-3},\ k=1.
\end{array}
\right.\label{HTo(2k)}\\
H(T_o(-2k,2))&=&\left\{
\begin{array}{l}
-z(a +a^{3}+\ldots +a^{2k-3})-(z+z^{-1})a^{2k-1} +z^{-1}a^{2k+1}, \ k>1,\\
-(z+z^{-1})a +z^{-1}a^{3}, \ k=1.
\end{array}
\right.\label{HTo(-2k)}
\end{eqnarray}
It follows that if $D$ is the connected sum of $\kappa^+$ positive torus links $T_o(2\alpha_j,2)$ and $\kappa^-$ negative torus links $T_o(-2\beta_i,2)$, then we have $r^+(D)=-\kappa^++\sum \alpha_j$ and $r^-(D)=-\kappa^-+\sum \beta_j$.
We have the following
\begin{eqnarray}
E(D)&=&s(D)-w(D)-1-2r^-(D)=\kappa^--\kappa^++2\sum\beta_i,\quad p_0^h(D)\in F,\label{TorusConnect3}\\
e(D)&=&-s(D)-w(D)+1+2r^+(D)=\kappa^+-\kappa^--2\sum\alpha_i,\quad  p_0^\ell(D)\in(-1)^{\kappa^++\kappa^-}F.
\label{TorusConnect4}
\end{eqnarray}

On the other hand, if $T_p(n, 2)$ and $T_p(-n,2)$ ($n\ge 2$) are the torus knots/links whose components (when $n$ is even) have parallel orientations and whose crossings are positive and negative respectively, then   
\begin{eqnarray}
H(T_p(n,2))&=&
f_{n+2} a^{1-n} - f_{n} a^{-1-n},
\label{HTp(n)}\\
H(T_p(-n,2))&=&
(-1)^{n-1}(f_{n+2} a^{n-1} - f_{n} a^{n+1}),
\label{HTp(-n)}
\end{eqnarray}
where $\{f_k\}$ is the Fibonacci like polynomial sequence defined by $f_{n+2}=zf_{n+1}+f_n$, $f_2=z^{-1}$ and $f_3=1$. We have $\deg(f_n)=n-3$. 
}
\end{remark}

\begin{remark}\label{mirror}{\em   It is well known that $H(\L^r,z,a)=H(\L,z,-a^{-1})$  for any link $\L$ where $\L^r$ is the mirror image of $\L$. From this one obtains that $E(\L^r)=-e(\L)$, $e(\L^r)=-E(\L)$, $p^h(\L^r)=(-1)^{e(\L)}p^\ell(\L)$ and $p^\ell(\L^r)=(-1)^{E(\L)}p^h(\L)$. It follows that $(E(\L)-e(\L))/2+1=(E(\L^r)-e(\L^r))/2+1$. Thus it suffices for us to prove our results in Theorem \ref{MT1} only for the cases when $\alpha_j>0$ for all $j$. If $\alpha_j<0$ for some $j$, then we must have $\beta_i>0$ for all $i$ by Definition \ref{TypeDef}, and we can obtain our results from $\L^r$ by switching the roles of $\kappa^+$, $\kappa^-$, $\alpha_j$ and $\beta_i$. Similarly, we shall only prove the cases $\kappa^+\ge \kappa^-$ for Type 2 pretzel links.
Since the braid index lower bounds are derived from the HOMFLY-PT polynomials, by Remark \ref{default}, we can re-arrange the order of the strips in the pretzel links to our liking. Thus for the next two sections, we shall further assume that $\alpha_1\ge \alpha_2 \ldots \ge \alpha_{\kappa^+}$ and $\beta_1\ge \beta_2 \ldots \ge \beta_{\kappa^-}$. In particular, $\alpha_{\kappa^+}=\min_{1\le j\le \kappa^+}\{\alpha_j\}$ and $\beta_{\kappa^-}=\min_{1\le i\le \kappa^-}\{\beta_i\}$.}
\end{remark}

\section{Braid index lower bounds for Type 1 pretzel links}\label{lowerbound1}

In this section, we prove that the formulas given in Theorem \ref{MT1}
are lower bounds of the braid indices of the corresponding
pretzel links. We do this by determining $E(D)$ and $e(D)$ for any
standard pretzel link diagram $D$ since every pretzel link is
represented by such a diagram. Our proofs are mostly induction based,
hence can be quite lengthy and tedious. In such a proof, it is
often the case that the beginning steps require a more detailed
analysis, with straight forward subsequent steps. In such cases, we
shall only provide the necessary details for the first step and write
``RLR" at the end of the first step, meaning ``the rest of the proof is
straight forward and is left to the reader". 

\medskip
Also,
some of the non alternating pretzel links, such as the ones with $\kappa^+=2$, $\kappa^-=1$ and $\beta_1=0$, in fact reduces to an alternating diagram that is no longer a pretzel link diagram. However, it is more convenient to include them in our consideration as this would make our induction based proof easier if we use them as our initial link diagrams in the induction process.

\medskip
For the sake of convenience, we shall use the notation $F$ ($-F$) to denote the
set of all nonzero Laurent polynomials of $z$ whose coefficients are
non-negative (non-positive). Notice that $F$ is closed under addition and multiplication. 

\medskip
\begin{proposition}\label{alter_1}
For $\L\in P_{1}(2\alpha_1+1,\ldots, 2\alpha_{\kappa^+}+1;0)$, we have
\begin{eqnarray*}
E(\L)&=&1-\kappa^+\,\,{\rm with} \ p^h(\L)= z^{\kappa^+-1},\\
e(\L)&=&-1-\kappa^+-2\sum_{j=1}^{\kappa^+}\alpha_j\,\,{\rm with} \,\, p^\ell(\L)=-f_{\kappa^+}.
\end{eqnarray*}
\end{proposition}

\begin{proof} Use induction on $\kappa^+$. For $\kappa^+=2$, $\L=T_o(2\alpha_1+2\alpha_2+2,2)$ and the statement follows from (\ref{HTo(2k)}). Assume that the statement holds for $2\le \kappa^+\le n$ for some $n\ge 2$ and consider the case $\kappa^+=n+1\ge 3$. If $\alpha_1 =\alpha_2 =\ldots =\alpha_{\kappa^+} =0$, then $\L=T_p(\kappa^+,2)$ and the statement holds by (\ref{HTp(n)}). 
Assume that the claim is true for some $\alpha_j \ge 0$, $1\ge j\ge \kappa^+$ and let us consider the case when one of the $\alpha_j$'s has been increased by one. Without loss of generality we may assume that $\alpha_1$ is increased from its current value $q$ to $q+1$. Use VP$^+$ at a  crossing in the first strip. The induction hypothesis applies to both $D_-$ and $D_0$, but It is easy to see that $E(\L)=-1+E(D_0)=1-\kappa^+$ which is contributed by $D_0$ only, and
$e(\L)=-1-\kappa^+-2\sum \alpha_j$ which is contributed by $D_-$ only, so we also have $p^h(\L)=zp^h(D_0)=z\cdot z^{\kappa^+-2}=z^{\kappa^+-1}$ and $p^\ell(\L)=p^\ell(D_-)=-f_{\kappa^+}$.
\end{proof}

\medskip
\begin{proposition}\label{T1}
If $\kappa^+=2$ and $\kappa^-=1$, then we have
$$
\begin{array}{lll}
&E(\L)=
\left\{
\begin{array}{ll}
-2, &\ {\rm if}\ \alpha_2>\beta_1,\\
-2\alpha_2+2\beta_1, &\ {\rm if}\ \alpha_2\le \beta_1,
\end{array}
\right.\quad &p^h(\L)\in F,\\
&e(\L)=
\left\{
\begin{array}{ll}
-2\alpha_1-2\alpha_2, &\ {\rm if}\ \beta_1=0,\\
-2-2\alpha_1-2\alpha_2, &\ {\rm if}\ \beta_1>0.\\
\end{array}
\right.
\quad &p^\ell(\L)\in \left\{
\begin{array}{ll}
-F, &\ {\rm if}\ \beta_1=0,\\
F, &\ {\rm if}\ \beta_1>0.\\
\end{array}
\right.
\end{array}
$$
\end{proposition}

\begin{proof} If $\beta_1=0$, then $\alpha_1\alpha_2>0$ and $\L\in P_1(2\alpha_1-1, 1, 2\alpha_2-1;0)= P_1(2(\alpha_1-1)+1, 1, 2(\alpha_2-1)+1;0)$ as shown in Figure \ref{beta_1=0}. The result then follows from Proposition \ref{alter_1}.

\begin{figure}[htb!]
\begin{center}
\includegraphics[scale=.45]{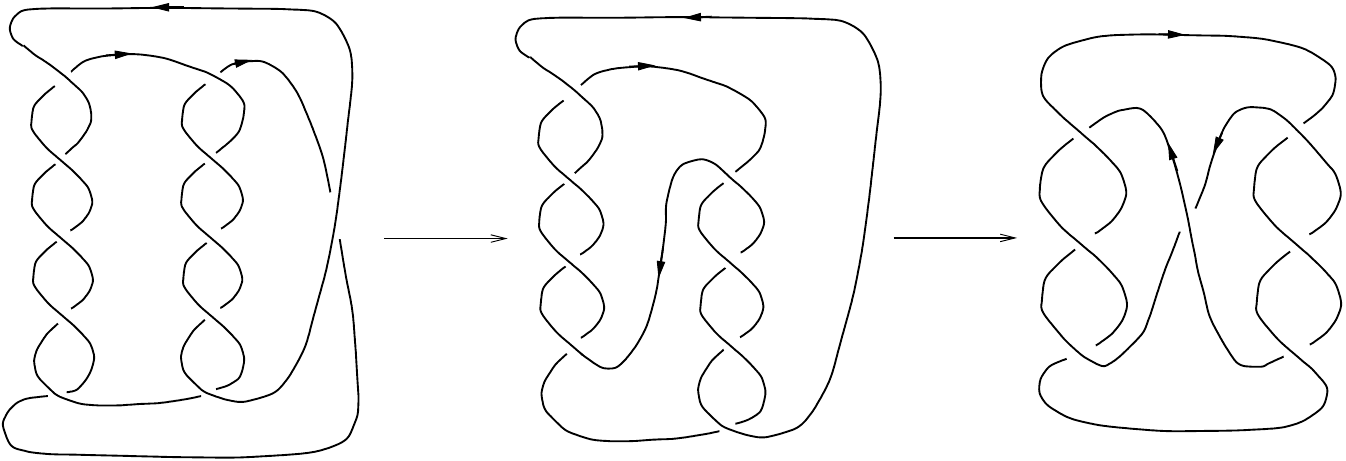} 
\end{center}
\caption{The non alternating pretzel link $P_1(5,5;-1)$ is equivalent to the alternating pretzel link diagram $P_1(3,3,1;0)$. }\label{beta_1=0}
\end{figure}

\medskip
Next we consider the case $\beta_1>0$ and $\alpha_2\le \beta_1$.  We claim that in fact 
$$
E(\L)=-2\alpha_2+2\beta_1,\quad e(\L)=-2-2\alpha_1-2\alpha_2
$$
with
$p^h(\L)\in F$ and $p^\ell(\L)\in F$.
Use induction on $\alpha_2$, starting with $\alpha_2=0$. If $\beta_1=1$, then we apply VP$^-$ on a negative crossing. $D_+$ is the unknot and $D_0$ is $T_o(2\alpha_1+2,2)$. $p^\ell(D_0)=-z^{-1}$ by (\ref{HTo(2k)}). We have 
\begin{eqnarray*}
E(a^{2}H(D_+))=2+E(D_+)=2,\ &&e(a^{2}H(D_+))=2+e(D_+)=2,\\
E(-zaH(D_0))=1+E(D_0)=0,\ &&e(-zaH(D_0))=1+e(D_0)=-2\alpha_1-2.
\end{eqnarray*}
Thus $E(\L)=2=2\beta_1$ with $p^h(\L)=1\in F$ and $e(\L)=-2-2\alpha_1$
 with $p^\ell(\L)=-zp^\ell(D_0)\in F$. Induction on $\beta_1$ then shows that $E(\L)=2\beta_1$ with $p^h(\L)=1\in F$ and $e(\L)=-2-2\alpha_1$ with $p^\ell(\L)=-zp^\ell(D_0)\in F$ for any $\beta_1\ge 1$ when $\alpha_2=0$. Assume now that $1\le \alpha_2\le \beta_1$ and that the statement is true for $\alpha_2-1$. Use VP$^+$ on a crossing in the $\alpha_2$ strip. The induction hypothesis applies to $D_-$ and $D_0$ is the torus link $T_0(2(\alpha_1-\beta_1),2)$, where the sign of $\alpha_1-\beta_1$ indicates the sign of the crossings in the torus link. In the special case of $\alpha_1-\beta_1=0$, $T_0(2(\alpha_1-\beta_1),2)$ is just an unlink with two components. (i) $\alpha_1>\beta_1$. By the induction hypothesis and direct calculations, we have
\begin{eqnarray*}
-2+E(D_-)&=&-2-2(\alpha_2-1)+2\beta_1=-2\alpha_2+2\beta_1,\\
-2+e(D_-)&=&-2-2-2\alpha_1-2(\alpha_2-1)=-2-2\alpha_1-2\alpha_2,\\
-1+E(D_0)&=&-2<0\le -2\alpha_2+2\beta_1,\\
-1+e(D_0)&=&-2-2(\alpha_1-\beta_1)> -2-2\alpha_2-2\alpha_1. 
\end{eqnarray*}
So in this case $E(\L)$ and $e(\L)$ are contributed by $D_-$ only. Hence 
$E(\L)=-2\alpha_2+2\beta_1$, $e(\L)=-2-2\alpha_1-2\alpha_2$ and
$
p^h(\L)=p^h(D_-)=1\in F$, $p^\ell(\L)=p^\ell(D_-)\in F
$
by the induction hypothesis. So the statement holds. (ii) $\alpha_1\le \beta_1$. The calculations involving $D_-$ remain unchanged, however the calculations involving  $D_0$ change since $D_0=T_0(2(\alpha_1-\beta_1),2)$ is either a trivial link with two components (when $\alpha_1=\beta_1$) or a torus link with negative crossings (when $\alpha_1<\beta_1$). By (\ref{HTo(-2k)}) we have 
\begin{eqnarray*}
-1+E(D_0)&=&2(\beta_1-\alpha_1)\le -2\alpha_2+2\beta_1, \ zp^h(D_0)=1\in F\\
-1+e(D_0)&=&0\,\,{\rm if}\, \alpha_1<\beta_1 \,\, {\rm or} \,\, -2 \, \,{\rm if}\, \alpha_1=\beta_1. 
\end{eqnarray*}
The first inequality follows because $\alpha_1\ge \alpha_2$ by our choice.
Thus  $e(\L)$ is contributed by $D_-$ only as in the case of $\alpha_1>\beta_1$, and $E(\L)=-2\alpha_2+2\beta_1$ even when $2(\beta_1-\alpha_1)= -2\alpha_2+2\beta_1$ (namely when $\alpha_1=\alpha_2$) since $zp^h(D_0)=1\in F$ and $p^h(D_-)=1\in F$.
This concludes the proof for the case $\beta_1>0$ and $\alpha_2\le \beta_1$. Notice that in the particular case $\alpha_2=\beta_1$, we have $E(\L)=0$.

\medskip
Finally we consider the case $\beta_1>0$ and $\alpha_1\ge \alpha_2\ge \beta_1$. We continue induction on $\alpha_2$, with the initial step $\alpha_2= \beta_1$ already established in the above. Assume that the statement is true for $\alpha_2-1\ge \beta_1$ and consider the case $\alpha_2$. Use VP$^+$ on a crossing in the $\alpha_2$ strip. Again, the induction hypothesis applies to $D_-$ and $D_0$ is the torus link $T_0(2(\alpha_1-\beta_1),2)$, where $\alpha_1-\beta_1>0$. We have
\begin{eqnarray*}
-2+E(D_-)&=&\left\{
\begin{array}{ll}
-2,\ &{\rm if}\ \alpha_2=\beta_1+1,\\
-4,\ &{\rm if}\ \alpha_2>\beta_1+1,\\
\end{array}
\right.\\
-2+e(D_-)&=&-2-2\alpha_1-2\alpha_2,\\
-1+E(D_0)&=&-2,\\
-1+e(D_0)&=&-2-2(\alpha_1-\beta_1)> -2-2\alpha_2-2\alpha_1. 
\end{eqnarray*}
So $e(\L)=-2+e(D_-)=-2-2\alpha_1-2\alpha_2$  with $p^\ell(\L)=p^\ell(D_-)\in F$. For the case $\alpha_2=\beta_1+1$, we have $-2+E(D_-)=-1+E(D_0)=-2$, however $p^h(D_-)\in F$ and $zp^h(D_0)=z^2\in F$ so $E(\L)=-2$ with $p^h(\L)=p^h(D_-)+zp^h(D_0)\in F$. For the case $\alpha_2>\beta_1+1$, $E(\L)=-1+E(D_0)=-2$ with $p^h(\L)=zp^h(D_0)\in F$. This concludes the proof of the proposition.
\end{proof}

\medskip
\begin{proposition}\label{T3}
If $\kappa^+\ge 3$, $\kappa^-=1$, then $e(\L)=-\kappa^+-2\sum \alpha_j$ and 
\begin{eqnarray*}
E(\L)&=&
\left\{
\begin{array}{ll}
2-\kappa^+-2\alpha_{\kappa^+}+2 \beta_1, & {\rm if}\ \alpha_{\kappa^+}\le \beta_1,\\
-\kappa^+, & {\rm if}\ \alpha_{\kappa^+}> \beta_1.
\end{array}
\right.
\end{eqnarray*}
Furthermore, $p^h(\L)\in F$, $p^\ell(\L)\in -F$ if $\beta_1=0$ and $p^\ell(\L)\in F$  if $\beta_1>0$.
\end{proposition}

\begin{proof} We shall use induction on $\kappa^+$, starting with $\kappa^+=3$.
First consider the case $\beta_1=0$. Starting at $\alpha_3=1$ and use VP$^+$ on a crossing in the $\alpha_3$ strip. $D_-=T_o(2(\alpha_1+\alpha_2+1),2)$ and Proposition \ref{T1} applies to $D_0$. By Remark \ref{toruslink_formula} and Proposition  \ref{T1}, we have the following:
\begin{eqnarray*}
-2+E(D_-)&=&-1+E(D_0)=-3,\ p^h(D_-)=z\in F,\ zp^h(D_0)\in F,\\
-2+e(D_-)&=&-5-2(\alpha_1+\alpha_2),\\
-1+e(D_0)&=&-1-2(\alpha_1+\alpha_2),\ zp^\ell(D_0)\in -F.
\end{eqnarray*}
Thus the statement holds. RLR for $\alpha_3\ge 1$ in general.

\medskip
Now consider the case $\beta_1>0$ and $\alpha_3=0$, starting with $\beta_1=1$. Use VP$^-$ on a negative crossing. $D_+=T_o(2(\alpha_1+\alpha_2+1),2)$ and $D_0\in P_1(2\alpha_1+1,2\alpha_1+1,1;0)$. By Proposition \ref{alter_1} and Remark \ref{toruslink_formula}, we have 
\begin{eqnarray*}
2+E(D_+)=1&>&-1=1+E(D_0),\\
2+e(D_+)=-1-2(\alpha_1+\alpha_2)&>&-3-2(\alpha_1+\alpha_2)=
1+e(D_0).
 \end{eqnarray*}
  Thus $E(\L)=1=2-\kappa^+-2\alpha_{\kappa^+}+2 \beta_1$, $p^h(\L)=p^h(D_+)=z\in F$, $e(\L)=-3-2(\alpha_1+\alpha_2)$ and $p^\ell(\L)=-zp^\ell(D_0)\in F$. RLR for $\beta_1>0$ in general.

\medskip
Next we consider the case $\beta_1>0$ and $\alpha_3>0$, starting with $\alpha_3=1$. Use VP$^+$ on a positive crossing in the $\alpha_3$ strip. The above case applies to $D_-$ and Proposition \ref{T1} applies $D_0$.  We have 
\begin{eqnarray*}
-2+E(D_-)&=&-2-1+2\beta_1=-1-2\alpha_3+2\beta_1,\\ 
-1+E(D_0)&=&\left\{
\begin{array}{ll}
-3,\ &{\rm if}\ \alpha_2>\beta_1,\\
-1-2\alpha_2+2\beta_1,\ &{\rm if}\ \alpha_2\le \beta_1,
\end{array}
\right.\\
-2+e(D_-)&=&-2-3-2(\alpha_1+\alpha_2)=-3-2\sum_{1\le j\le 3} \alpha_j,\\
-1+e(D_0)&=&-3-2(\alpha_1+\alpha_2). 
\end{eqnarray*}
Thus $e(\L)=-3-2\sum_{1\le j\le 3} \alpha_j$ with $p^\ell(\L)=p^\ell(D_-)\in F$. Furthermore,
$E(\L)=-3+2\beta_1$ with $p^h(\L)=p^h(D_-)+zp^h(D_0)\in F$ if $\alpha_2=\alpha_3=1$. If $\alpha_2>\alpha_3=1$, then $-1-2\alpha_2+2\beta_1<-1-2\alpha_3+2\beta_1$. Also, $-3<-1-2\alpha_3+2\beta_1=-3+2\beta_1$. Thus in this case $E(\L)=-3+2\beta_1=-1-2\alpha_3+2\beta_1$ with $p^h(\L)=p^h(D_-)\in F$. So the statement of the theorem holds.

\medskip
Assume now that the statement holds for $\alpha_3=n$ where $1\le n\le \beta_1-1$, we need to show that the statement still holds for $\alpha_3=n+1$. Again use VP$^+$ on a positive crossing in the $\alpha_3$ strip. The induction hypothesis applies to $D_-$ and  $D_0$ is as in the case of $\alpha_3=1$. We have 
\begin{eqnarray*}
-2+e(D_-)&=&-2-3-2(\alpha_1+\alpha_2+n)=-3-2\sum_{1\le j\le 3} \alpha_j\\
&<&-3-2(\alpha_1+\alpha_2)=-1+e(D_0),
\end{eqnarray*}
so $e(\L)=-3-2\sum \alpha_j$ with $p^\ell(\L)=p^\ell(D_-)\in F$ by the induction hypothesis.
On the other hand, by comparing
$
-2+E(D_-)=-3-2n+2\beta_1=-1-2\alpha_3+2\beta_1
$
with $-1+e(D_0)=-1-2\alpha_2+2\beta_1$, we see that $E(\L)=-1-2\alpha_3+2\beta_1$ with $p^h(\L)\in F$ since $\alpha_3\le \alpha_2$ and $p^h(D_-)\in F$, $zp^h(D_0)\in F$.

\medskip
We have now shown that the statement is true for any $\alpha_3\le \beta_1$. If we continue induction on $\alpha_3$ for $\alpha_3\ge\beta_1$, then the above discussion has already established the initial step  $\alpha_3=\beta_1$. RLR. This completes the proof for the initial step $\kappa^+=3$. 

\medskip
Assume that the statement is true for $\kappa^+\ge n\ge 3$, we now need to prove that the statement holds for $\kappa^+=n+1$. The proof is very similar to the case of $\kappa^+=3$ and is left to the reader. 
\end{proof}

\medskip
\begin{proposition}\label{T5}
If $\kappa^+\ge 2$, $\kappa^-\ge 2$, $\alpha_j>0$ for $1\le j\le \kappa^+$, then
\begin{eqnarray}
E(\L)&=&1-(\kappa^+-\kappa^-)+2\sum \beta_i,\label{T5e1}\\
e(\L)&=&
\left\{
\begin{array}{ll}
1-(\kappa^+-\kappa^-)-2\sum \alpha_j,\ {\rm if}\ \tau=\kappa^+-1,\\
-1-(\kappa^+-\kappa^-)-2\sum \alpha_j,\ {\rm otherwise},
\end{array}
\right.
\label{T5e2}
\end{eqnarray}
where $\tau$ is the number of $\beta_i$'s that equal to zero.
\end{proposition}

\begin{proof} We shall prove (\ref{T5e1}) first, together with the claim that $p^h(\L)\in (-1)^{\kappa^-}F$. Use induction on $\gamma=\kappa^--\tau$, which is the number of positive $\beta_i$'s. Let us first establish the case $\gamma=0$. Use induction on $\kappa^-\ge 2$, starting at $\kappa^-=2$. Use VP$^-$ on a negative crossing. Proposition \ref{alter_1} applies to $D_+$ and Proposition \ref{T3} applies to $D_0$. 
Thus $2+E(D_+)=2+1-\kappa^+=1-(\kappa^+-\kappa^-)>1+E(D_0)=1-\kappa^+$.
Therefore $E(\L)=1-(\kappa^+-\kappa^-)$, $p^h(\L)=p^h(D_+)\in F=(-1)^{\kappa^-}F$, and the statement holds for any $\kappa^+$. Assume that the statement holds for $2\le \kappa^-\le q$ for some $q\ge 2$, consider the case $\kappa^-=q+1$. Again apply VP$^-$ to a negative crossing. The induction hypothesis applies to $D_0$, and also to $D_+$ if $q\ge 3$. If $q=2$, then Proposition \ref{T3} applies to $D_+$. If $q\ge 3$, then by the induction hypothesis we have
$2+E(D_+)=1+E(D_0)=1-(\kappa^+-\kappa^-)$. The result follows since $p^h(D_+)\in (-1)^{\kappa^--2}F=(-1)^{\kappa^-}F$ and $-zp^h(D_0)\in -(-1)^{\kappa^--1}F=(-1)^{\kappa^-}F$. If $q=2$, then by Proposition \ref{T3} we have $2+E(D_+)=2-\kappa^+<1-(\kappa^+-\kappa^-)=4-\kappa^+=1+E(D_0)$, so the result also holds. This proves (\ref{T5e1}) for $\gamma=0$.

\medskip
 Assume that (\ref{T5e1}) holds for $\gamma\le n$ for some $n\ge 0$ and consider the case $\gamma=n+1\ge 1$. Use induction on $\beta_\gamma$, starting with $\beta_\gamma=1$. Use VP$^-$ on a negative crossing in the $\beta_{\gamma}$ strip. While the induction hypothesis applies to $D_+$, it only applies to $D_0$ if $\kappa^-\ge 3$. In the case that $\kappa^-=2$, either Proposition \ref{T1} ($\kappa^+=2$) or Proposition \ref{T3} ($\kappa^+\ge 3$) applies to $D_0$. Thus for $\kappa^-\ge 3$, we have
$$
2+E(D_+)=1-(\kappa^+-\kappa^-)+2\sum\beta_i>-1-(\kappa^+-\kappa^-)+2\sum\beta_i=1+E(D_0),
$$
hence $E(\L)=2+E(D_+)=1-(\kappa^+-\kappa^-)+2\sum\beta_i$ with $p^h(\L)=p^h(D_+)\in (-1)^{\kappa^-}F$. If $\kappa^-=2$, then $E(D_0)\le -\kappa^++2\beta_1$ by Propositions \ref{T1} and \ref{T3}. It follows that 
$$
2+E(D_+)=1-(\kappa^+-\kappa^-)+2\sum\beta_i\ge 3-\kappa^++2\beta_1>1-\kappa^++2\beta_1\ge 1+E(D_0).
$$
Induction on $\beta_{\gamma}$ then extends this result to any $\beta_{\gamma}\ge 1$ trivially, as the calculations about $D_0$ remains the same. This completes the induction and (\ref{T5e1}) is proved.

\medskip
Now we will provide a detailed proof for (\ref{T5e2}), by proving several claims that will lead to the proof of (\ref{T5e2}).

\medskip\noindent
Claim 1. If $\gamma=0$, then 
\begin{eqnarray}
e(\L)&=&
\left\{
\begin{array}{ll}
-1-(\kappa^+-\kappa^-)-2\sum \alpha_j,\ {\rm if}\ \kappa^-\not=\kappa^+-1,\\
-2\sum \alpha_j,\ {\rm if}\ \kappa^-=\kappa^+-1,
\end{array}
\right.\label{Claim1e1}
\end{eqnarray}
with
\begin{eqnarray}
p^\ell(\L)&=&\left\{
\begin{array}{ll}
-f_{\kappa^+-\kappa^-},\ {\rm if}\ \kappa^-\le \kappa^+-2,\\
(-1)^{\kappa^++\kappa^-+1}f_{\kappa^--\kappa^++2}\ {\rm if}\ \kappa^-\ge\kappa^+,\\
-\kappa^-=-(\kappa^+-1)\ {\rm if}\ \kappa^-=\kappa^+-1.
\end{array}
\right.\label{Claim1e2}
\end{eqnarray}
\begin{proof} Consider first the case $1\le \kappa^-\le \kappa^+-2$ (so it is necessary that $\kappa^+\ge 3$). Use induction on $\kappa^-$ and VP$^-$ to a negative crossing at each step, start with $\kappa^-=1$. By Proposition \ref{alter_1} we have
$2+e(D_+)=1+e(D_0)=-\kappa^+-2\sum \alpha_j$, 
$
p^\ell(D_+)=-f_{\kappa^++1}
$
and 
$-zp^\ell(D_0)=zf_{\kappa^+}$.
Thus $e(\L)=-\kappa^+-2\sum \alpha_j=-1-(\kappa^+-\kappa^-)-2\sum \alpha_j$ with $p^\ell(\L)=-(f_{\kappa^++1}-zf_{\kappa^+})=-f_{\kappa^+-1}=-f_{\kappa^+-\kappa^-}$.
For $\kappa^-=2$, we have (by Proposition \ref{alter_1} and the above case of $\kappa^-=1$)
\begin{eqnarray*}
2+e(D_+)&=&2-1-\kappa^+-2\sum \alpha_j=1-\kappa^+-2\sum \alpha_j,\quad p^\ell(D_+)=-f_{\kappa^+},\\
1+e(D_0)&=&1-\kappa^+-2\sum \alpha_j,\quad -zp^\ell(D_0)=zf_{\kappa^+-1}.
\end{eqnarray*}
It follows that 
 $e(\L)=1-\kappa^+-2\sum \alpha_j=-1-(\kappa^+-\kappa^-)-2\sum \alpha_j$ with $p^\ell(\L)=-f_{\kappa^+}+zf_{\kappa^+-1}=-f_{\kappa^+-2}=-f_{\kappa^+-\kappa^-}$. Assume that the statement is true for $\kappa^-=n\ge 2$, consider the case $\kappa^-=n+1$. This time the induction hypothesis applies to both $D_0$ and $D_+$. A similar argument then leads to $e(\L)=-1-(\kappa^+-\kappa^-)-2\sum \alpha_j$ with $p^\ell(\L)=-f_{\kappa^+-\kappa^-+2}+zf_{\kappa^+-\kappa^-+1}=-f_{\kappa^+-\kappa^-}$. 
  
 \medskip
Next we consider the cases $\kappa^-=\kappa^+-1$ and $\kappa^-=\kappa^+$ together. For $\kappa^+=2$ and $\kappa^-=1$, the statement follows from Proposition \ref{alter_1} by observing that $P_1(2\alpha_1+1,2\alpha_2+1;-1)=P_1(2\alpha_1-1,2\alpha_2-1,1;0)$, see Figure \ref{beta_1=0}. For the case $\kappa^+=\kappa^-=2$, apply VP$^-$ to a negative crossing. $D_+=T_o(2\alpha_1+2\alpha_2+2,2)$ and $D_0$ is the case $\kappa^+=2$ and $\kappa^-=1$. Thus $2+e(D_+)=-1-2(\alpha_1+\alpha_2)<1+e(D_0)=1-2(\alpha_1+\alpha_2)$, hence $e(\L)=-1-2(\alpha_1+\alpha_2)$ with $p^\ell(\L)=p^\ell(D_+)=-f_2=-z^{-1}$. This establishes the first step of the induction. In general, assume that the statement is true for $\kappa^+=n\ge 2$, $\kappa^+-1\le \kappa^-\le \kappa^+$, let us consider the case $\kappa^+=n+1$, starting with $\kappa^-=\kappa^+-1=n$ and $\alpha_{n+1}=1$. Use 
 VP$^+$ on a positive crossing in the $\alpha_{n+1}$ strip. $D_-$ has $n$ positive strips and $n-1$ negative strips and $D_0$ has $n$ positive and $n$ negative strips. By the induction hypothesis we have 
 $$
-2+ e(D_-)=-1+e(D_0)=-2-2\sum_{1\le j\le n} \alpha_j=-2\sum_{1\le j\le n+1} \alpha_j,
 $$
 and
 $$
 p^\ell(\L)=p^\ell(D_-)+zp^\ell(D_0)=-(n-1)+z(-z^{-1})=-n=-\kappa^-.
 $$
 RLR for any $\alpha_{n+1}\ge 1$. This proves the case of $\kappa^-=n$. 
For $\kappa^-=\kappa^+=n+1$, use VP$^-$ on a negative crossing. The first part of the proof above applies to $D_+$ (which has $n+2$ positive strips and $n$ negative strips) and the induction hypothesis applies to $D_0$ (which has $n+1$ positive strips and $n$ negative strips). We have
\begin{eqnarray*}
 &&2+e(D_+)= 2-1-(n+2-n)-2\sum \alpha_j=-1-2\sum \alpha_j\\
 &<&
 1+e(D_0)=1-1-(n+1-n)-2\sum \alpha_j= 1-2\sum \alpha_j.
\end{eqnarray*}
 It follows that $e(\L)=-1-2\sum \alpha_j$ with $p^\ell(\L)=p^\ell(D_+)=-f_{n+2-n}=-f_2=-z^{-1}$ by the first part of the proof above. This concludes the proof of the cases $\kappa^-=\kappa^+-1$ and $\kappa^-=\kappa^+$.
 
\medskip
Finally we consider the cases $\kappa^-\ge \kappa^+$. Use induction on $\kappa^-$. The initial step $\kappa^-= \kappa^+$ is already established above. For $\kappa^-=\kappa^++1$, apply VP$^-$ on a negative crossing. By the above results, we have 
$1+e(D_0)=-2\sum \alpha_j < 2+e(D_+)=2-2\sum \alpha_j$, thus $e(\L)=-2\sum \alpha_j=-1-(\kappa^+-\kappa^-)-2\sum \alpha_j$ with 
$$p^\ell(\L)=-zp^\ell(D_0)=-z(-z^{-1})=1=f_3=(-1)^{\kappa^+-\kappa^-+1}f_{\kappa^--\kappa^++2},$$
as desired. Assume that the statement holds for $\kappa^-\le n$ for some $n\ge \kappa^++1$, consider the case $\kappa^-= n+1\ge \kappa^++2$. Apply VP$^-$ on a negative crossing. The induction hypothesis now applies to both $D_0$ and $D_-$. We have 
 \begin{eqnarray*}
&&1+e(D_0)=1-1-(\kappa^+-n)-2\sum \alpha_j= -1-(\kappa^+-\kappa^-)-2\sum \alpha_j\\
&=& 2+e(D_+)=2-1-(\kappa^+-(n-1))-2\sum \alpha_j= -1-(\kappa^+-\kappa^-)-2\sum \alpha_j.
\end{eqnarray*}
Thus $e(\L)=-1-(\kappa^+-\kappa^-)-2\sum \alpha_j$ with 
 \begin{eqnarray*}
p^\ell(\L)&=&p^\ell(D_+)-zp^\ell(D_0)=(-1)^{\kappa^++(n-1)+1}f_{n-1-\kappa^++2}-z(-1)^{\kappa^++n+1}f_{n-\kappa^++2}\\
&=&(-1)^{\kappa^++n}(f_{n+1-\kappa^+}+zf_{n+2-\kappa^+})=(-1)^{\kappa^++\kappa^-+1}f_{\kappa^--\kappa^++2}.
\end{eqnarray*}
\end{proof}

 \medskip\noindent
Claim 2. If $1\le \kappa^-\le \kappa^+-1$ and $\gamma\ge 1$, then (\ref{Claim1e1}) still holds, but
$
p^\ell(\L)=(-1)^{\gamma+1}z^{\gamma} f_{\kappa^+-\tau}.
$
\begin{proof} Use induction on $\gamma$. The initial step $\gamma=0$ is guaranteed by Claim 1 if $1\le \kappa^-\le \kappa^+-2$. Assume that the statement holds for $\gamma=n$ for some $n\ge 0$, consider the case $\gamma=n+1$.
Start with $\beta_{n+1}=1$. Apply VP$^-$ on a negative crossing in the $\beta_{n+1}$ strip. By the induction hypothesis (which applies to both $D_+$ and $D_0$), we have 
\begin{eqnarray*}
&&1+e(D_0)=1-1-(\kappa^+-n)-2\sum \alpha_j=-1-(\kappa^+-\kappa^-)-2\sum \alpha_j\\
&<&2+e(D_+)=2-1-(\kappa^+-(n+1))-2\sum \alpha_j=1-(\kappa^+-\kappa^-)-2\sum \alpha_j.
\end{eqnarray*} 
Thus $e(\L)=-1-(\kappa^+-\kappa^-)-2\sum \alpha_j$ with 
$$
p^\ell(\L)=-zp^\ell(D_0)=-z(-1)^{(\gamma-1)+1}z^{\gamma-1} f_{\kappa^+-\tau}=(-1)^{\gamma+1}z^{\gamma} f_{\kappa^+-\tau}.
$$
RLR for $\beta_{n+1}\ge 1$ in general. This proves the case $1\le \kappa^-\le \kappa^+-2$. In the case that $ \kappa^-= \kappa^+-1$, we have to start with $\gamma=1$ (so $\tau=\kappa^--1=\kappa^+-2$). For $\beta_1=1$, use VP$^-$ on a negative crossing in the $\beta_1$ strip. By Claim 1, we see that 
$$
e(\L)=1+e(D_0)=1-1-2-2\sum \alpha_j=-2-2\sum \alpha_j,
$$
with $p^\ell(\L)=-zp^\ell(D_0)=-z(-f_{\kappa^+-(\kappa^+-2)})=zf_2=(-1)^{\gamma+1}z^\gamma f_{\kappa^+-\tau}$. RLR for $\beta_1\ge 1$ in general. 
Now assume the statement is true for $1\le \gamma\le q$ for some $q\ge 1$ and consider the case of $\gamma=q+1$ with $\beta_{q+1}=1$. Use VP$^-$ on a negative crossing in the $\beta_{q+1}$ strip. The inductive assumption applies to $D_+$ and the first part of the proof above applies to $D_0$ (which has $\kappa^+-2$ negative strips but the same $\tau$ value as that of $D$), we have 
$$
2+e(D_+)=-2\sum \alpha_j>
1+e(D_0)= 1-1-(\kappa^+-(\kappa^+-2))-2\sum \alpha_j=-2-2\sum \alpha_j.
$$
 Thus
$e(\L)=1+e(D_0)=-2-2\sum \alpha_j$ with
$$
p^\ell(\L)=-zp^\ell(D_0)=-z(-1)^{q+1}z^q f_{\kappa^+-\tau}=(-1)^{q+2}z^{q+1} f_{\kappa^+-\tau}=(-1)^{\gamma+1}z^{\gamma} f_{\kappa^+-\tau},
$$
as desired. RLR for $\beta_{q+1}\ge 1$ in general. This completes the proof of the claim.
\end{proof}

\medskip\noindent
Claim 3. If $\kappa^-= \kappa^+$, then
\begin{eqnarray}
e(\L)&=&
\left\{
\begin{array}{ll}
-1-(\kappa^+-\kappa^-)-2\sum \alpha_j=-1-2\sum \alpha_j,\ {\rm if}\ \tau<\kappa^+-1,\\
1-(\kappa^+-\kappa^-)-2\sum \alpha_j=1-2\sum \alpha_j,\ {\rm if}\ \tau=\kappa^+-1,
\end{array}
\right.\label{Claim3e1}\\
p^\ell(\L)&=&\left\{
\begin{array}{ll}
(-z)^{\gamma-1}, \ {\rm if}\  \tau<\kappa^+-1,\\
(\kappa^+-1)z-z^{-1},\ {\rm if}\ \tau=\kappa^+-1, \beta_1=1,\\
(\kappa^+-1)z,\ {\rm if}\  \tau=\kappa^+-1, \beta_1>1.
\end{array}
\right.\label{Claim3e2}
\end{eqnarray}
\begin{proof}
Start with $\gamma=1$ and $\beta_1=1$. Use VP$^-$ on a negative crossing in the $\beta_1$ strip. By Claim 1, we have $e(\L)=2+e(D_+)=1+e(D_0)=1-2\sum \alpha_j$ with 
$$
p^\ell(\L)=p^\ell(D_+)-zp^\ell(D_0)=-z^{-1}-z(-\kappa^++1)=(\kappa^+-1)z-z^{-1}.
$$
For $\beta_1\ge 2$, it is easy to see that $e(\L)$ is always given by $1+e(D_0)=1-2\sum \alpha_j$ hence $p^\ell(\L)=-zp^\ell(D_0)=(\kappa^+-1)z$. 

\medskip
For $\gamma=2$, again starts with $\beta_2=1$ and use VP$^-$ on a negative crossing in the $\beta_2$ strip. The above case of $\gamma=1$ applies to $D_+$ and Claim 2 applies to $D_0$. Thus we have 
$$
2+e(D_+)=3-2\sum \alpha_j>1+e(D_0)=1-2-2\sum \alpha_j=-1-2\sum \alpha_j.$$
It follows that 
$e(\L)=-1-2\sum \alpha_j$ with $p^\ell(\L)=-zp^\ell(D_0)=-z=(-z)^{\tau-1}$. RLR for $\beta_2\ge 2$ and for the general case $\gamma\ge 2$ by the same argument.
\end{proof}

\medskip\noindent
Claim 4. If $\tau\ge \kappa^+$, then
\begin{eqnarray}
e(\L)&=&
-1-(\kappa^+-\kappa^-)-2\sum \alpha_j, \quad p^\ell(\L)=(-1)^{\kappa^++\kappa^-+1}z^{\gamma}f_{\tau-\kappa^++2}. \label{Claim4e1}
\end{eqnarray}
\begin{proof} Use induction on $\gamma=\kappa^--\tau$. Note that the initial step $\gamma=0$ has been established in Claim 1. Assume that the statement holds for $\gamma=q$ for some $q\ge 0$, consider the case $\gamma=q+1$. 
Start with $\beta_{q+1}=1$ and apply VP$^-$ to a negative crossing in the $\beta_{q+1}$ strip. The induction hypothesis applies to both $D_+$ and $D_0$. We have
\begin{eqnarray*}
&&2+e(D_+)=1-(\kappa^+-\kappa^-)-2\sum \alpha_j\\
&>&
1+e(D_0)=1-1-(\kappa^+-(\kappa^--1))-2\sum \alpha_j=-1-(\kappa^+-\kappa^-)-2\sum \alpha_j,
\end{eqnarray*}
thus $e(\L)=-1-(\kappa^+-\kappa^-)-2\sum \alpha_j$ with 
$$
p^\ell(\L)=-zp^\ell(D_0)=-z(-1)^{\kappa^++(\kappa^--1)+1}z^{q}f_{\tau-\kappa^++2}=(-1)^{\kappa^++\kappa^-+1}z^{q+1}f_{\tau-\kappa^++2}.
$$
RLR for $\beta_{q+1}\ge 2$. 
\end{proof}

\medskip\noindent
Claim 5. If $\tau= \kappa^+-1$, then
\begin{eqnarray}
e(\L)&=&
1-(\kappa^+-\kappa^-)-2\sum \alpha_j,\quad  p^\ell_0(\L)=(-1)^{\gamma+1}(\kappa^+-1)z^{\gamma}. \label{Claim5e1}
\end{eqnarray}
\begin{proof} 
Use induction on $\gamma$. The statement holds for $\gamma=0$ by Claim 1. Assume that it holds for $\gamma=q$ for some $q\ge 0$. Then for $\gamma=q+1$, start with $\beta_\gamma=1$ and apply VP$^-$ to a negative crossing in the $\beta_\gamma$ strip. The induction hypothesis applies to $D_0$ since it still has $\tau=\kappa^+-1$ and Claim 4 applies to $D_+$ as it has a $\tau$ value of $\kappa^+$. We have
\begin{eqnarray*}
&&2+e(D_+)=2-1-(\kappa^+-\kappa^-)-2\sum \alpha_j=1-(\kappa^+-\kappa^-)-2\sum \alpha_j\\
&=&1+e(D_0)=1-1-(\kappa^+-(\kappa^--1))-2\sum \alpha_j.
\end{eqnarray*}
By Claim 4 and the induction hypothesis, we have 
\begin{eqnarray*}
p^\ell(D_+)&=&(-1)^{\kappa^+-\kappa^-+1}z^{q}f_{\kappa^+-\kappa^++2}=(-1)^{\gamma}z^{\gamma-2},\\
-zp_0^\ell(D_0)&=&-z(-1)^{(\gamma-1)+1}(\kappa^+-1)z^{\gamma-1}=(-1)^{\gamma+1}(\kappa^+-1)z^{\gamma}.
\end{eqnarray*}
The statement follows. RLR for $\beta_\gamma\ge 2$. 
\end{proof}

\medskip\noindent
Claim 6. If $\tau< \kappa^+-1$ and $\kappa^-\ge \kappa^+$, then 
\begin{eqnarray*}
e(\L)&=&
-1-(\kappa^+-\kappa^-)-2\sum \alpha_j,\quad  p^\ell(\L)=(-z)^{\gamma-1}. 
\end{eqnarray*}
\begin{proof} 
Use induction on $\kappa^-$. We note that the initial step $\kappa^-=\kappa^+$ is guaranteed by Claim 3. Assume that the statement holds for $\kappa^-=q$ for some $q\ge \kappa^+$, consider the case $\kappa^-=q+1$, starting with $\beta_{\gamma}=1$. and apply VP$^-$ to a negative crossing in the $\beta_\gamma$ strip. The induction hypothesis applies to $D_0$, and either the induction hypothesis (in the case that $\tau<\kappa^+-2$) or Claim 5 (in the case that $\tau=\kappa^+-2$) applies to $D_+$. But either way we have
\begin{eqnarray*}
&&2+e(D_+)\ge 2-1-(\kappa^+-\kappa^-)-2\sum \alpha_j=1-(\kappa^+-\kappa^-)-2\sum \alpha_j\\
&>&1+e(D_0)=1-1-(\kappa^+-(\kappa^--1))-2\sum \alpha_j=-1-(\kappa^+-\kappa^-)-2\sum \alpha_j.
\end{eqnarray*}
Thus $e(\L)=-1-(\kappa^+-\kappa^-)-2\sum \alpha_j$ with $p^\ell(\L)=-zp^\ell(D_0)=-z(-z)^{(\gamma-1)-1}=(-z)^{\gamma-1}$. RLR for $\beta_{\gamma}\ge 2$.
\end{proof}

\medskip
Combining Claims 1 to 6, which have covered all possible cases for $\kappa^+\ge 2$ and $\kappa^-\ge 2$, completes the proof of Proposition \ref{T5}.
 \end{proof}

\medskip
\begin{theorem}\label{Type1LowboundTheorem}
The formulas given in Theorem \ref{MT1} are lower bounds of the braid indices of the corresponding pretzel links.
\end{theorem}

\begin{proof} This is established by computing $(E(\L)-e(\L))/2+1$ using the results obtained in this section. In the following, we list the results from which the lower bound for each formula in Theorem \ref{MT1} is derived from, but leave the calculations of $(E(\L)-e(\L))/2+1$ to the reader.

\medskip
\begin{tabular}{ll}
Formula (\ref{MT1e1}):& Proposition \ref{T1},\\
Formula (\ref{MT1e2}):& Propositions \ref{T1}, \ref{T3},\\
Formula (\ref{MT1e3}):& Propositions \ref{T1}, \ref{T3},\\
Formula (\ref{MT1e4}):& Proposition \ref{T5},\\
Formula (\ref{MT1e5}):& Propositions \ref{alter_1}, \ref{T5}.
\end{tabular}
\end{proof}

\section{Braid index lower bounds for Type 2 pretzel links}\label{lowerbound2}

\begin{proposition}\label{T4.1}
If $\kappa^+=2$ and $\kappa^-=1$, then
\begin{eqnarray*}
E(\L)&=&
\left\{\begin{array}{ll}
-2, &{\rm if}\ \alpha_2> \beta_1,\\
2(\beta_1-\alpha_2),\ &{\rm if}\ \alpha_2\le \beta_1,
\end{array}
\right.\\
e(\L)&=&-2(\alpha_1+\alpha_2).
\end{eqnarray*}
\end{proposition}

\begin{proof} We shall prove the statement with the claim that $p^h(\L)\in F$ and $p^\ell(\L)\in F$. Consider first the case $\alpha_2\le \beta_1$.  Start with $\alpha_2=1$ and apply VP$^+$ to a crossing in the $\alpha_2$ strip. $D_0=T_o(2\alpha_1-2\beta_1,2)$ and $D_-=T_o(2\alpha_1,2)\#T_o(-2\beta_1,2)$. The statement follows from direct calculation (which we shall leave to the reader). Assume that the statement is true for $\alpha_2\le n$ for some $n$ such that $1\le n \le \beta_1-1$, consider the case $\alpha_2=n+1\le \beta_1$. Again apply VP$^+$ to a crossing in the $\alpha_2$ strip. The induction hypothesis applies to $D_-$, while $D_0$ remains $T_o(2\alpha_1-2\beta_1,2)$. By the induction hypothesis $-2+E(D_-)=-2+2(\beta_1-n)=2(\beta_1-\alpha_2)$, and direction calculation shows that $-1+E(D_0)\le 2(\beta_1-\alpha_2)$. Since $p^h(D_-)\in F$ (by the induction hypothesis) and $zp^h(D_0)\in F$, we have $E(\L)=2(\beta_1-\alpha_2)$ with $p^h(\L)\in F$ as desired. On the other hand, $-2+e(D_-)=-2-2(\alpha_1+n)=-2(\alpha_1+\alpha_2)$, and direction calculation shows that $-1+E(D_0)> -2(\alpha_1+\alpha_2)$. Thus $e(\L)=-2(\alpha_1+\alpha_2)$ with $p^\ell(\L)\in F$. This completes the proof for the case $\alpha_2\le \beta_1$. 

\medskip
Next we consider the case $\alpha_2\ge \beta_1$ using induction on $\alpha_2$. Notice that the above has established the initial step $\alpha_2= \beta_1$.
Assume that the statement is true for $\alpha_2\le n$ for some $n$ such that $n \ge \beta_1$, consider the case $\alpha_2=n+1> \beta_1$. Again apply VP$^+$ to a crossing in the $\alpha_2$ strip. Direction calculation shows that $-1+E(D_0)=-2$ with $zp^h(D_0)=z^2\in F$ since $\alpha_1>\beta_1$. By the induction hypothesis $-2+E(D_-)=-2$ if $\beta_1=n$ or $-4$ if $\beta_1<n$. Either way we have $E(\L)=-2$ with $p^h(\L)\in F$ since $p^h(D_-)\in F$ (by the induction hypothesis). $e(\L)=-2(\alpha_1+\alpha_2)$ with $p^\ell(\L)\in F$ is shown as before. This completes the proof for the case $\alpha_2\ge \beta_1$. 
\end{proof}

\medskip
We now generalize Proposition \ref{T4.1} to Proposition \ref{T4.2} below.

\medskip
\begin{proposition}\label{T4.2}
For $\kappa^+\ge 2$ and $\kappa^-=1$, we have
\begin{eqnarray*}
E(\L)&=&\left\{
\begin{array}{ll}
-\kappa^+\ &{\rm if}\ \alpha_{\kappa^+}> \beta_1,\\
2-\kappa^+-2\alpha_{\kappa^+}+2\beta_1\ &{\rm if}\ \alpha_{\kappa^+}\le \beta_1,
\end{array}
\right.\\
e(\L)&=&2-\kappa^+-2\sum \alpha_j
\end{eqnarray*}
with $p^h(\L)\in F$ and $p^\ell(\L)\in (-1)^{\kappa^+}F$.
\end{proposition}

\begin{proof}
We shall use double induction on $\kappa^+$ and $\alpha_{\kappa^+}$. The initial step for $\kappa^+=2$ is proved in Proposition \ref{T4.1}. Assume that the statement of the theorem holds for $\kappa^+\le n$ for some $n\ge 2$, let us now consider the case $\kappa^+=n+1$. Start with $\alpha_{n+1}=1$ and apply VP$^+$ to a crossing in the $\alpha_{n+1}$ strip. The induction hypothesis applies to $D_0$ (with $\kappa^+=n$) and $D_-$ is the connected sum of  $T_o(2\alpha_1,2)$, ..., $T_o(2\alpha_n,2)$ and $T_o(-2\beta_1,2)$. We have (by Remark \ref{toruslink_formula})
\begin{eqnarray*}
-2+E(D_-)&=&-2-n+2\beta_1+1=2-(n+1)+2(\beta_1-\alpha_{n+1}),\\
-1+E(D_0)&=&
\left\{
\begin{array}{ll}
-(n+1),\ &{\rm if}\ \alpha_n> \beta_1,\\
2-(n+1)+2(\beta_1-\alpha_n),\ &{\rm if}\ \alpha_n\le \beta_1,
\end{array}
\right.\\
-2+e(D_-)&=&-2+1-n-2\sum_{1\le j\le n}\alpha_j=2-\kappa^+-2\sum_{1\le j\le \kappa^+}\alpha_j,\\
-1+e(D_0)&=&-1+2-n-2\sum_{1\le j\le n}\alpha_j=4-\kappa^+-2\sum_{1\le j\le \kappa^+}\alpha_j.
\end{eqnarray*}
By comparison we see that the statement of the theorem holds with $p^h(\L)=p^h(D_-)\in F$ if $\alpha_n>1(=\alpha_{n+1})$ and $p^h(\L)=p^h(D_-)+zp^h(D_0)\in F$ if $\alpha_n=1$, and $p^\ell(\L)=p^\ell(D_-)\in (-1)^{\kappa^+}F$. 
This proves the case of $\alpha_{\kappa^+}=1$. It is straight forward to see that the above formula for $e(\L)$ holds in general for all $\alpha_{\kappa^+}$ and $(-1)^{\kappa^+}p^\ell(\L)\in F$. The case for $E(\L)$ is a little more complicated and we need to consider the cases $\alpha_{\kappa^+}\le \beta_1$  and $\alpha_{\kappa^+}>\beta_1$. 
Notice that in the above case we have $\alpha_{\kappa^+}=1\le \beta_1$ automatically. Assume that the statement is true for $\alpha_{\kappa^+}\ge q\ge 1$ with $q<\beta_1$, let us consider the case $\alpha_{\kappa^+}=q+1$. Apply VP$^+$ to a crossing in the $\alpha_{\kappa^+}$ strip. $D_0$ is the same as before and the hypothesis applies to $D_-$. In this case we have
\begin{eqnarray*}
-2+E(D_-)&=&-\kappa^+-2q+2\beta_1>-\kappa^+,\\
-1+E(D_0)&=&
\left\{
\begin{array}{ll}
2-\kappa^+-2\alpha_n+2\beta_1,\ &{\rm if}\ \alpha_n\le \beta_1,\\
-\kappa^+,\ &{\rm if}\ \alpha_n> \beta_1.
\end{array}
\right.
\end{eqnarray*}
Since $\alpha_n\ge \alpha_{\kappa^+}=q+1$, $2-\kappa^+-2\alpha_n+2\beta_1\le -\kappa^+-2q+2\beta_1$. Thus $E(\L)=-\kappa^+-2q+2\beta_1=2-\kappa^+-2\alpha_{\kappa^+}+2\alpha_1$ with $p^h(\L)\in F$ since $p^h(D_-)\in F$ and $zp^h(D_0)\in F$. This shows that the statement is true for all $\alpha_{\kappa^+}$ such that $\alpha_{\kappa^+}-1<\beta_1$. In particular, the statement is also true for $\alpha_{\kappa^+}=\beta_1$. Now let us consider the case $\alpha_{\kappa^+}=\beta_1+1$. This time we have $-2+E(D_-)=-\kappa^+$ and $-1+E(D_0)=-\kappa^+$ since $\alpha_n\ge \alpha_{\kappa^+}>\beta_1$. It follows that $E(\L)=-\kappa^+$ with $p^h(\L)=p^h(D_-)+zp^h(D_0)\in F$. From here it is easy to see that $E(\L)=-\kappa^+$ with $p^h(\L)=zp^h(D_0)\in F$ for any  $\alpha_{\kappa^+}>\beta_1+1$. \end{proof}

\medskip
\begin{proposition}\label{T4.3}
For $\kappa^+\ge 2$ and $\kappa^-\ge 2$, we have
\begin{eqnarray*}
E(\L)&=&-1+\kappa^--\kappa^++2\sum \beta_i,\ p^h(\L)\in F,\\
e(\L)&=&1+\kappa^--\kappa^+-2\sum \alpha_j,\ (-1)^{\kappa^++\kappa^-+1}p^\ell(\L)\in F.
\end{eqnarray*}
\end{proposition}

\begin{proof} We use induction on $\kappa^-$. The case $\kappa^-=1$ has been proven by Propositions \ref{T4.1} and \ref{T4.2}, so the base case starts with $\kappa^-=2$. Start with $\beta_2=1$ and apply $VP^-$ to the $\beta_2$ strip. Using (\ref{TorusConnect3}) and (\ref{TorusConnect4}) in Remark \ref{toruslink_formula} for $D_+$ and Propositions \ref{T4.1} and \ref{T4.2} for $D_0$, we have
\begin{eqnarray*}
2+E(D_+)&=&2+1-\kappa^++2\beta_1=-1+\kappa^--\kappa^++2\sum_{1\le j\le \kappa^-}\beta_i>1+E(D_0),\\
1+E(D_0)&=&\left\{
\begin{array}{ll}
1-\kappa^+=-1+\kappa^--\kappa^+\ &{\rm if}\ \alpha_{\kappa^+}> \beta_1,\\
3-\kappa^+-2\alpha_{\kappa^+}+2\beta_1=1+\kappa^--\kappa^+-2\alpha_{\kappa^+}+2\beta_1\ &{\rm if}\ \alpha_{\kappa^+}\le \beta_1,
\end{array}
\right.\\
2+e(D_+)&=&2-1+\kappa^++2\sum_{1\le j\le \kappa^+}\alpha_j=3-\kappa^-+\kappa^++2\sum_{1\le j\le \kappa^+}\alpha_j,\\
1+e(D_0)&=&3-\kappa^+-2\sum_{1\le j\le \kappa^+}\alpha_j=1+\kappa^--\kappa^++2\sum_{1\le j\le \kappa^+}\alpha_j.
\end{eqnarray*}
The statement follows since $p^h(D_+)\in F$ and $-zp^\ell(D_0)\in -(-1)^{\kappa^+}F=(-1)^{\kappa^++\kappa^-+1}F$. Now assume that the statement holds for $\beta_2\le q$ for some $q\ge 1$, then for $\beta_2= q+1$, the inductive statement applies to $D_+$ while the powers for $D_0$ remain unchanged.
\begin{eqnarray*}
2+E(D_+)&=&-1+\kappa^--\kappa^++2\sum_{1\le j\le \kappa^-}\beta_i,\\
2+e(D_+)&=&3+\kappa^--\kappa^+-2\sum_{1\le j\le \kappa^+} \alpha_j.
\end{eqnarray*}
$E(\L)=-1+\kappa^--\kappa^++2\sum_{1\le j\le \kappa^-}\beta_i$ is contributed by $D_+$ only with $p^h(\L)=p^h(D_+)\in F$, and $e(\L)=1+\kappa^--\kappa^++2\sum_{1\le j\le \kappa^+}\alpha_j$ is contributed by $D_0$ only with $p^\ell(\L)=-zp^h(D_0)\in (-1)^{\kappa^++\kappa^-+1}F$. 
Assume that the statement is true for $\kappa^-\le n$ for some $n\ge 1$ and consider $\kappa^-=n+1$. Start with $\beta_{n+1}=1$ and apply $VP^-$ to the $\beta_{n+1}$ strip. Using (\ref{TorusConnect3}) and (\ref{TorusConnect4}) in Remark \ref{toruslink_formula} for $D_+$ and the induction hypothesis for $D_0$, we have
\begin{eqnarray*}
2+E(D_+)&=&-1+\kappa^--\kappa^++2\sum_{1\le j\le \kappa^-}\beta_i>-3+\kappa^--\kappa^++2\sum_{1\le j\le \kappa^-}\beta_i=1+E(D_0),\\
2+e(D_+)&=&3+\kappa^+-\kappa^--2\sum_{1\le j\le \kappa^+}\alpha_j>1+\kappa^--\kappa^+-2\sum_{1\le j\le \kappa^+} \alpha_j=1+e(D_0).
\end{eqnarray*}
The inequality in the second equation follows from the fact that we assume $\kappa^+\ge \kappa^-$, see Remark \ref{mirror}.
The statement follows since $p^\ell(\L)=-zp^\ell(D_0)\in -(-1)^{1+\kappa^++n}F=(-1)^{\kappa^++\kappa^-+1}F$ and $p^h(\L)=p^\ell(D_-)\in F$. Now assume that the statement holds for $\beta_{n+1}\le q$ for some $q\ge 1$, then for $\beta_{n+1}= q+1$, the inductive statement applies to $D_+$ while the powers for $D_0$ remain unchanged.
$E(\L)=-1+\kappa^--\kappa^++2\sum_{1\le j\le \kappa^-}\beta_i$ is contributed by $D_+$ only with $p^h(\L)=p^h(D_+)\in F$, and $e(\L)=1+\kappa^--\kappa^++2\sum_{1\le j\le \kappa^+}\alpha_j$ is contributed by $D_0$ only with $p^\ell(\L)=-zp^h(D_0)\in (-1)^{\kappa^++\kappa^-+1}F$. 
\end{proof}

\medskip
\begin{theorem}\label{Type2LowboundTheorem}
The formulas given in Theorem \ref{MT2} are lower bounds of the braid indices of the corresponding pretzel links.
\end{theorem}

\begin{proof} Again, this is established by computing $(E(\L)-e(\L))/2+1$ using the results obtained in this section. In the following, we list the results from which the lower bound for each formula in Theorem \ref{MT1} is derived from and leave the calculations to the reader.

\medskip
\begin{tabular}{ll}
Formula (\ref{MT2e1}):& The proof of the Type M2 case in \cite[Theorem 4.7]{Diao2021}\\
Formula (\ref{MT2e2}):& Proposition \ref{T4.2}\\
Formula (\ref{MT2e3}):& Proposition \ref{T4.2}\\
Formula (\ref{MT2e4}):& Proposition \ref{T4.3}
\end{tabular}
\end{proof}

\section{preparations for establishing the upper bounds}\label{upper_prep}

\begin{figure}[!htb]
\includegraphics[scale=0.8]{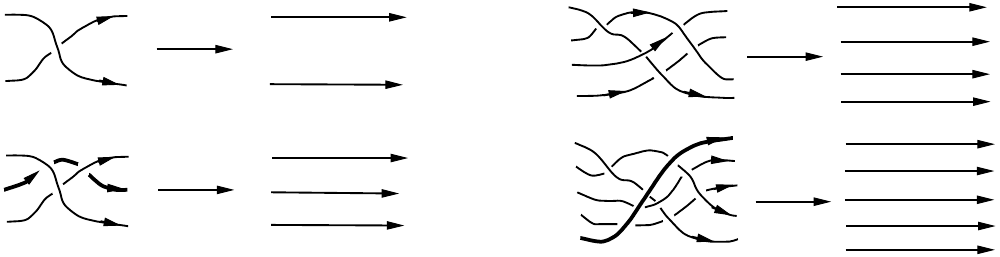}
\caption{Examples of a strand passing through a local braid with the same orientation. Top: the original local braid and its Seifert circle decomposition; Bottom: the local braid after a strand has been re-routed through it and the Seifert circle decomposition after that.}\label{local_braid}
\end{figure}

\medskip
\begin{remark}\label{open_braid}{\em 
In order to establish an upper bound for the braid index of a pretzel link, we use a result of Yamada~\cite{Ya}, which states that the number of Seifert circles in a link diagram $D$ is an upper bound of the braid index of the link. In our approach, we will start from a standard diagram of a pretzel link, and apply a sequence of ambient isotopy moves to reduce the number of Seifert circles in the diagram, such that the number of Seifert circles at the end of the process matches the formulas given in Theorems \ref{MT1} and \ref{MT2}. These ambient isotopy moves are usually in the form of re-routing an over strand or an under strand at a crossing, and the following observation will help the reader to understand the Seifert circle decomposition after such a re-routing is performed. As demonstrated in Figure \ref{local_braid}, if our diagram has a local braid structure and we re-route a strand over (or under) this local braid in the same direction of the local braid, then we will end up with another local braid with one more string.}
\end{remark}

\begin{figure}[!htb]
\includegraphics[scale=0.5]{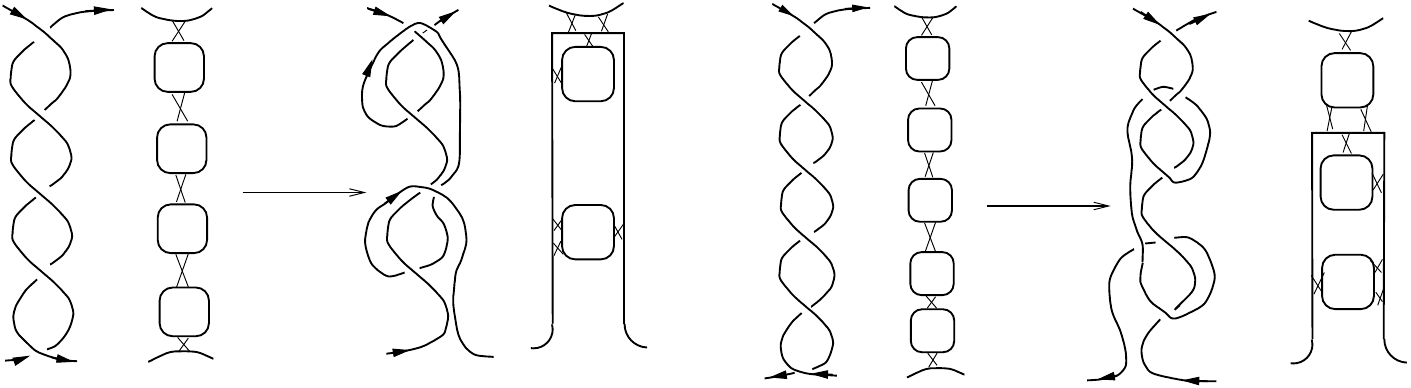}
\caption{The multiple MP moves performed on a string of Seifert circles connected by lone crossings. Left: The number of Seifert circles is reduced from $2k$ to $k$ (with $k=2$ here); Right: The number of Seifert circles is reduced from $2k+1$ to $k+1$ (with $k=2$ here).}\label{MP_move}
\end{figure}

\medskip
There are two special types of moves that will be used repeatedly so we will define them here. The first move, called the {\em Murasugi-Przytycki move} (MP move for short)~\cite{Diao2021}, reduces the number of Seifert circles in a diagram $D$ by one at a {\em lone crossing} (a crossing that is the only crossing between two Seifert circles). See Figure \ref{MP_move} for an illustration of this move. Notice that the re-routed strand in the MP move always enters a local braid that with the same orientation as the local braid. If a diagram has a lone crossing, then we can always reduce its number of Seifert circles by one by performing a MP move at the lone crossing. In our case, we will be applying multiple MP moves to strings of Seifert circles created by the vertical strips in the Type 1 and Type 2 pretzel links. In the case of a Type 1 pretzel link, each strip contributes an even number of Seifert circles. Say this number is $2k$ with $k\ge 1$, then we can reduce it by $k$ as shown in Figure \ref{MP_move}. On the other hand, a strip in a Type 2 pretzel link contributes an odd number of Seifert circles. Say this number is $2k+1$ with $k\ge 0$, then we can also reduce it by $k$ as shown in Figure \ref{MP_move}.

\begin{remark}\label{repeated_MP_move}{\em
For a Type 1 pretzel link $P_1(2\alpha_1+1,\ldots,2\alpha_{\kappa^+}+1;-(2\beta_1+1),\ldots,-(2\beta_{\kappa^-}+1)$, it is easy to see that it has $2+2\sum\alpha_j+2\sum\beta_i$ Seifert circles and $\sum\alpha_j+\sum\beta_i$ MP moves, hence it is always possible to represent $\L$ by a diagram with $2+\sum\alpha_j+\sum\beta_i$ Seifert circles. For a Type 2 pretzel link $P_2(2\alpha_1,\ldots,2\alpha_{\kappa^+};-(2\beta_1),\ldots,-(2\beta_{\kappa^-})$, it is easy to see that it has $2+2\sum\alpha_j+2\sum\beta_i-(\kappa^++\kappa^-)$ Seifert circles and $\sum\alpha_j+\sum\beta_i-(\kappa^++\kappa^-)$ MP moves, leaving us with  a diagram with $2+\sum\alpha_j+\sum\beta_i$ Seifert circles. However this diagram still contains at least one lone crossings, and so after one additional move hence it is always possible to represent $\L$ by a diagram with $1+\sum\alpha_j+\sum\beta_i$ Seifert circles (see \cite{Diao2021}).
}
\end{remark}

\begin{figure}[!htb]
\includegraphics[scale=0.8]{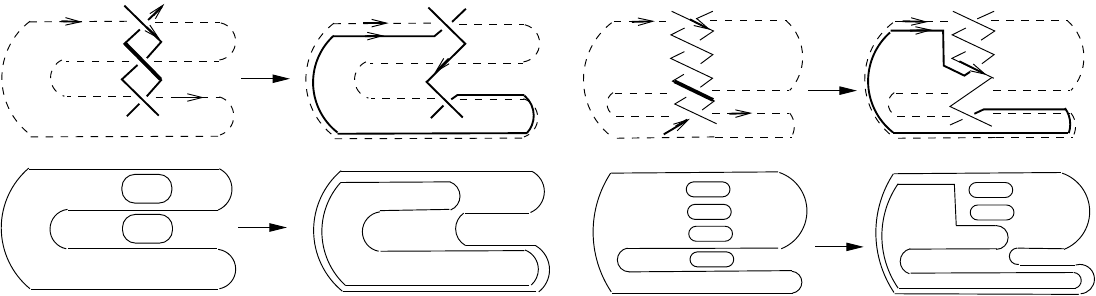}
\caption{Top: Eligible crossings where a special move can be made. The over strands to be re-routed are highlighted, the dashed line indicates the Seifert circle $C$ going under the vertical strip of crossings; Bottom: The Seifert circle decompositions before and after a special move. }\label{SMove}
\end{figure}

The second move, a new move introduced in this paper which we shall call a {\em special move} (S
move for short), applies at a crossing in a strip of $2k+1\ge 3$ crossings
as shown in Figure \ref{SMove}. In this case, the strands of a Seifert circle $C$ (shown by dashed lines in Figure \ref{SMove}) go underneath the crossings in the strip, and the rerouting of the over strand at the second crossing in the strip, as shown in Figure \ref{SMove}, is guided by the Seifert circle $C$.  A special move also reduces the number of Seifert circles in the diagram by one.

\section{Braid index upper bounds for Type 1 pretzel links}\label{upperbound1}

In this section, we will show that, by direct construction, the formulas given in Theorems \ref{MT1} are the upper bounds of the corresponding pretzel links. We first make the following observation.

\begin{figure}[htb!]
\includegraphics[scale=0.9]{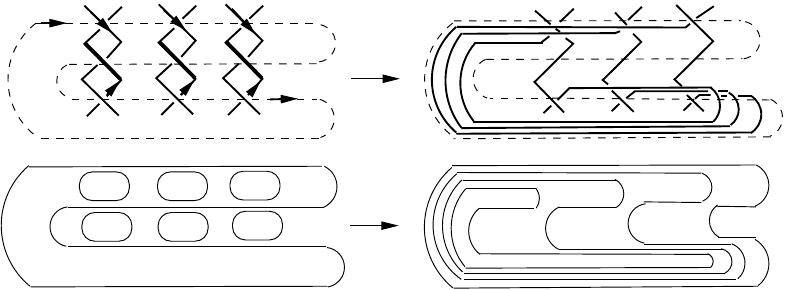}
\caption{Shown is an example of compatible multiple special moves. }
\label{repeat_longmoves}
\end{figure}

\begin{figure}[htb!]
\includegraphics[scale=0.8]{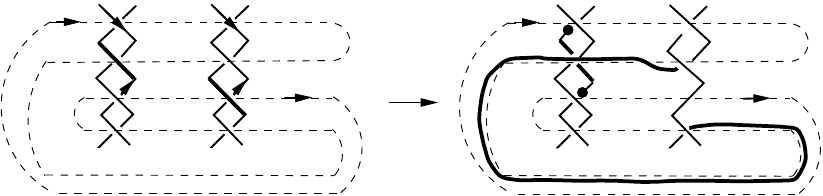}
\caption{An S move can prevent another S move as shown. Notice that in this case the rerouted (over) strand has to over cross the second over strand (shown by thickened curves).}
\label{non_repeatable}
\end{figure}

\begin{remark}\label{repeat_long_move_remark}{\em In the case that multiple S moves exist at different crossings which all involve the same long Seifert circle, a rerouted strand may prevent another S move. For example, if both moves involve over strands and one of the rerouted strand has to go over the other over strand, then the second over strand becomes unusable as demonstrated by Figure \ref{non_repeatable}. We say that two or more S moves are {\em compatible} if the rerouted strand of any of these moves will not prevent the other S moves. Figure \ref{repeat_longmoves} shows an example of multiple S moves that are compatible. 
}
\end{remark}

\medskip
\begin{remark}\label{permutation}{\em In general, changing the order of the tangles appearing in a Montesinos link may change the link type. However, the link type of a Montesinos link will not change if we perform a cyclic permutation of the tangles, or if we reverse the order of the tangles \cite[Section 12D]{Burde}. Furthermore, if a strip contains a single crossing, then we can place the strip in any position by using the flype moves, which will not change the link type.}
\end{remark}

In the following, we shall provide concrete constructions showing that the expressions given in formulas (\ref{MT1e1}) to (\ref{MT1e5}) are also upper bounds for the braid indices of the corresponding pretzel links. Our constructions will only use the properties mentioned in Remark \ref{permutation}, and will not depend on the order of the positive or negative strips in the pretzel links.

\medskip
Proof of (\ref{MT1e1}): We have $\L\in P_1(2\alpha_1+1, 2\alpha_2+1;-1)= P_1(2(\alpha_1-1)+1, 1, 2(\alpha_2-1)+1;0)$ as shown in Figure \ref{beta_1=0}. We have $\b(\L)=2+(\alpha_1-1)+(\alpha_2-1)=\alpha_1+\alpha_2$ by \cite[Theorem 4.7]{Diao2021}. This proves the case of  (\ref{MT1e1}).

\begin{figure}[htb!]
\includegraphics[scale=.6]{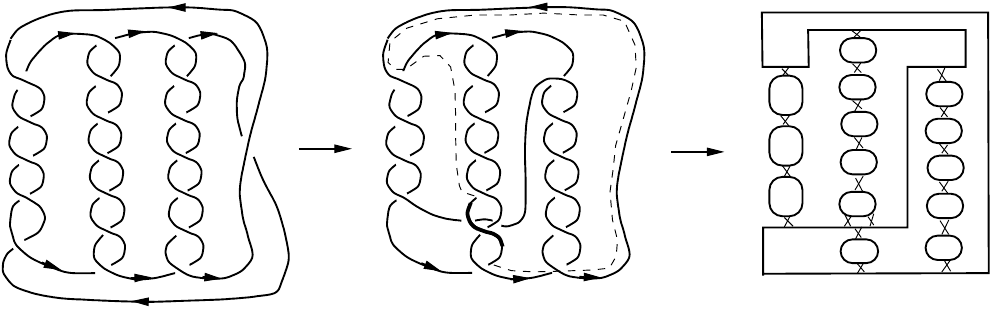}
\caption{The move taking the diagram $D$ of $\L\in P_1(5,7,7;-1)$ to $\tilde{D}$ with $s(\tilde{D})=-1+2\sum\alpha_j=15$ and $-2+\sum\alpha_j=6$ S/MP moves. The over strand at which an S move can be made is highlighted with the intended re-routing of the strand indicated by dashed line.}
\label{short_move}
\end{figure}

\medskip
Proof for upper bound in  (\ref{MT1e2}): Consider first the case $\kappa^+\ge 3$ and $\beta_1=0$. Let $D$ be a standard diagram of $\L$. Use flypes if necessary, we can assume that the only negative crossing is at the right side of the pretzel link as shown in Figure \ref{short_move}. The move shown in Figure \ref{short_move} reduces one Seifert circle each from the $\alpha_1$ strip and the $\alpha_{\kappa^+}$ strip, and combines the two large Seifert circles into one. Notice that there now $2\alpha_1-1$ Seifert circles from the $\alpha_1$ strip and $2\alpha_{\kappa^+}-1$ Seifert circles from the $\alpha_{\kappa^+}$ strip, with reduction numbers of $\alpha_1-1$ and $\alpha_{\kappa^+}-1$ respectively. For $2\le j\le \kappa^+-1$, each $\alpha_j$ strip still contributes $2\alpha_j$ Seifert circles, but each strip has one S move and $\alpha_j-1$ MP moves. By Remark \ref{repeat_long_move_remark}, these S moves are repeatable. Thus we can obtain a diagram of $\L$ with $-1+2\sum \alpha_j-(-2+\sum \alpha_j)=1+\sum \alpha_j$. This settles the case of $\beta_1=0$.

\begin{figure}[htb!]
\includegraphics[scale=.6]{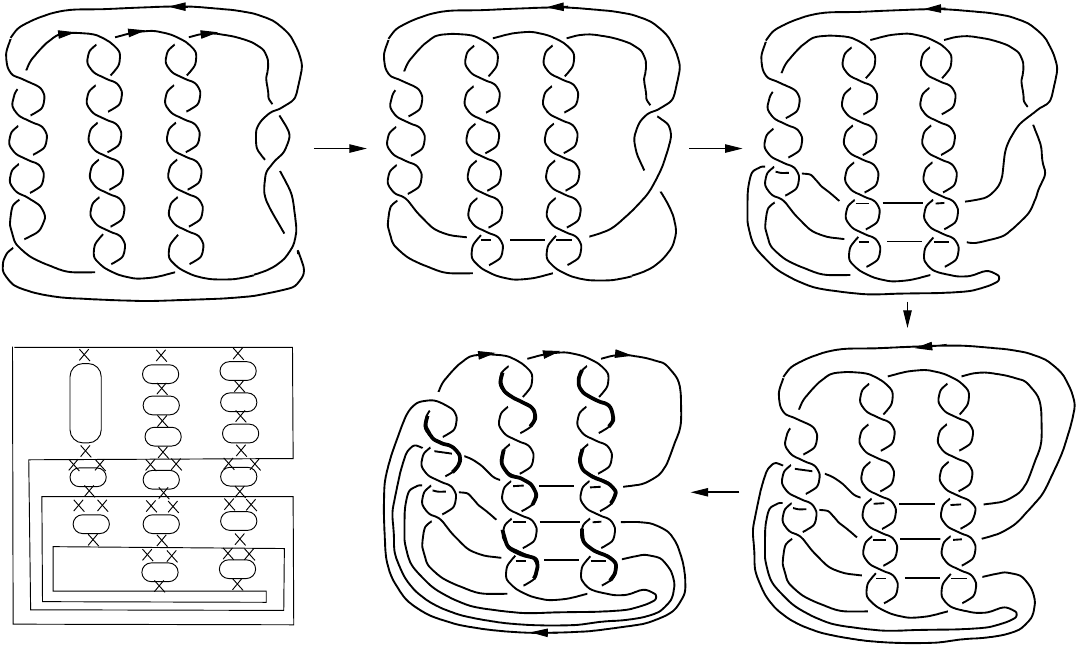}
\caption{The long moves taking the diagram $D$ of $\L\in P_1(5,7,7;-3)$ to $\tilde{D}$ with $s(\tilde{D})=2\sum \alpha_j=16$ and $1+\sum \alpha_j=7$ S/MP moves (2 MP moves and 5 compatible S moves). The over strands at which these moves can be made are highlighted by thick lines.}
\label{beta1small}
\end{figure}

Now consider the case $\kappa^-=1$ with $0<\beta_1<\min_{1\le j\le \kappa^+}\{\alpha_j\}$. The construction here works for the case of $\kappa^+=2$ as well. In this case, a standard diagram $D$ of $\L$ can be deformed to a new diagram $\tilde{D}$ by a sequence of re-routing moves (which we shall call {\em long} moves or L moves for short) as shown in Figure \ref{beta1small}. The condition $\beta_1<\min_{1\le j\le \kappa^+}\{\alpha_j\}$ allows the re-routed strands to intersect each local braid structure with  the same orientations. The Seifert circle decomposition of $\tilde{D}$ contains a long Seifert circle as shown in Figure \ref{beta1small}. 
$\tilde{D}$ retains the Seifert circles contributed by the $\alpha_j$ strips for $j\ge 2$, and $2\alpha_1-1$ of the $2\alpha_1$ Seifert circles contributed by the $\alpha_1$ strip. That is, $s(\tilde{D})=2\sum \alpha_j$. Furthermore, there are 
$\alpha_1-1$ S/MP moves in the $\alpha_1$ strip, and $\alpha_j$ S/MP moves in the $\alpha_j$ strip for $1< j\le \kappa^+$. The top strands at these crossings are marked by thickened curves in Figure \ref{beta1small}. One can verify that these S moves are compatible hence they are all realizable. Thus we can reduce the number of Seifert circles in $\tilde{D}$ by $-1+\sum \alpha_j$. Hence $\L$ can be represented by a closed braid with $2\sum \alpha_j-(-1+\sum\alpha_j)=1+\sum\alpha_j$ strings. This completes the proof for the upper bound in the case of (\ref{MT1e2}). Note that there is no S move for bottom crossing of the $\alpha_1$ strip, as this is not compatible with the other S-moves, see Figure \ref{non_repeatable}.

\begin{figure}[htb!]
\includegraphics[scale=.6]{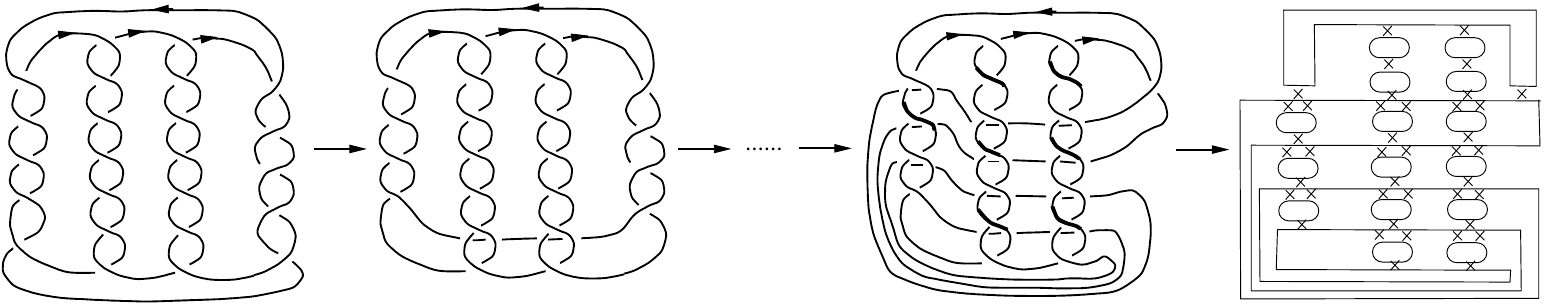}
\caption{The long moves taking the diagram $D$ of $\L\in P_1(5,7,7;-5)$ to $\tilde{D}$ with $s(\tilde{D})=17$ and 7 additional (compatible) S/MP moves.}
\label{beta1equal}
\end{figure}

 \medskip
Proof for upper bound in (\ref{MT1e3}): 
In the case that $\beta_1=\min_{1\le j\le \kappa^+}\{\alpha_j\}$, as shown in Figure \ref{beta1equal}, the long moves change $D$ to a diagram $\tilde{D}$ with a long Seifert circle such that $s(\tilde{D})=1+2\sum \alpha_j$. $\tilde{D}$ retains all Seifert circles contributed by each $\alpha_j$ strip for $2\le j\le \kappa^+$, and all but one Seifert circles contributed by the $\alpha_1$ strip. For $2\le j\le \kappa^+$, each $\alpha_j$ strip has $\alpha_j$ S/MP moves ($\beta_1$ S moves and $\alpha_j-\beta_1$ MP moves to be precise), but the $\alpha_1$ strip has only $\alpha_1-1$ S/MP moves ($\beta_1-1$ S moves and $\alpha_1-\beta_1$ MP moves, see Remark \ref{repeat_long_move_remark}). Thus $\tilde{D}$ has a total of 
$-1+\sum \alpha_j$ (compatible) S/MP moves, as shown in Figure \ref{beta1equal}. Hence $\L$ has a diagram with 
 $1+2\sum \alpha_j-(-1+\sum \alpha_j)=2+\sum \alpha_j$ Seifert circles.
 
\medskip
Similarly, in the case that $\beta_1>\min_{1\le j\le \kappa^+}\{\alpha_j\}$, say $\alpha_{\kappa^+}=\min_{1\le j\le \kappa^+}\{\alpha_j\}$, then the long moves will remove $2\alpha_{\kappa^+}$ crossings from the $\beta_1$ strip and reduce the number of Seifert circles by $1+2\alpha_{\kappa^+}$. There are $-1+\sum \alpha_j+(\beta_1-\alpha_{\kappa^+})$ additional 
 S/MP moves available in the resulting diagram. Thus we will have a representative diagram of $\L$ with 
 $2+\sum \alpha_j+\beta_1-\alpha_{\kappa^+}=2+\sum \alpha_j+\beta_1-\min_{1\le j\le \kappa^+}\{\alpha_j\}$ Seifert circles. This proves the upper bound for the case of (\ref{MT1e3}).
 
\begin{figure}[htb!]
\includegraphics[scale=.6]{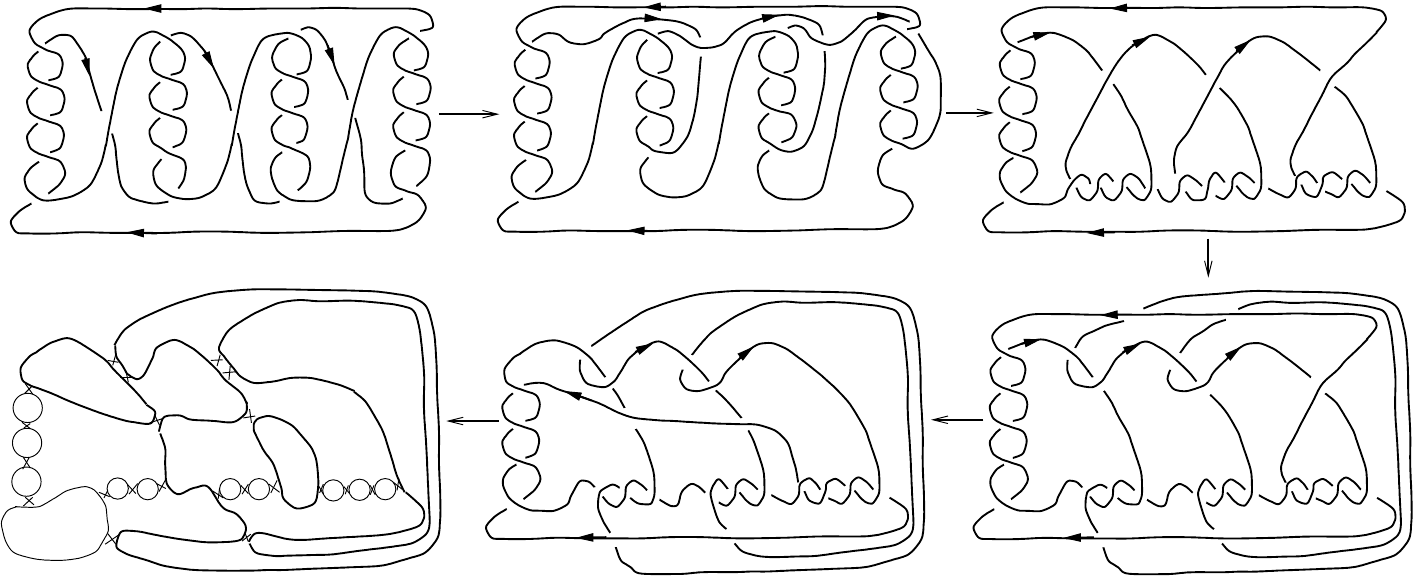}
\caption{The sequence of moves taking a diagram $D$ of $\L\in P_1(5,5,5,5;-1,-1,-1)$ to a diagram $\tilde{D}$ such that $s(\tilde{D})=3-\kappa^++2\sum \alpha_j=15$. Notice that there are $2-\kappa^++\sum \alpha_j=6$ MP moves available hence $\L$ has a diagram with $9=1+\sum \alpha_j$ Seifert circles.}
\label{extra_reduction}
\end{figure}

\medskip
Proof for upper bound in (\ref{MT1e4}): We only need to consider the case $\tau=\kappa^+-1$, that is, the number of $\beta_i$'s that are equal to zero is $\kappa^+-1$. Consider first the case $\kappa^-=\tau=\kappa^+-1$. Using flype moves if necessary, we can choose a representative diagram $D$ of $\L$ as shown in the starting diagram of Figure \ref{extra_reduction}. A sequence of moves shown in Figure \ref{extra_reduction} will result in a diagram $\tilde{D}$ with $s(\tilde{D})=3-\kappa^++2\sum \alpha_j$ Seifert circles as follows:
First,  for $2\le j\le \kappa^+$, each $\alpha_j$ strip is converted into a horizontal strip containing $2\alpha_j$ crossings (top row of Figure \ref{extra_reduction}). Next, for $2\le j\le \kappa^+-1$, each $\alpha_j$ strip is converted into a horizontal strip containing $2\alpha_j-1$ crossings contributing $2\alpha_j-2$ Seifert circles to the diagram, while the $\alpha_1$ and $\alpha_{\kappa^+}$ strips  contributes a sequence of $2\alpha_1-1$ and $2\alpha_{\kappa^+}-1$ Seifert circles to $\tilde{D}$ (bottom row of Figure \ref{extra_reduction}).
Furthermore, $\tilde{D}$ contains $\kappa^+-2$ concentric large Seifert circles with three additional Seifert circles not concentric to these large Seifert circles as demonstrated in Figure \ref{extra_reduction}. 
Furthermore, there are still $\alpha_j-1$ MP moves available in each $\alpha_j$ strip for $2\le j\le \kappa^+-1$, and there are $\alpha_1$ and $\alpha_{\kappa^+}$ MP moves each for the $\alpha_1$ and $\alpha_{\kappa^+}$ strips respectively. Thus the total number of MP moves in $\tilde{D}$ is $2-\kappa^++\sum \alpha_j$. It follows that $\L$ can be represented by a diagram with $1+\sum \alpha_j$ Seifert circles. This proves the upper bound for (\ref{MT1e4}) for the case $\kappa^-=\tau=\kappa^+-1$.
 
\begin{figure}[htb!]
\includegraphics[scale=.5]{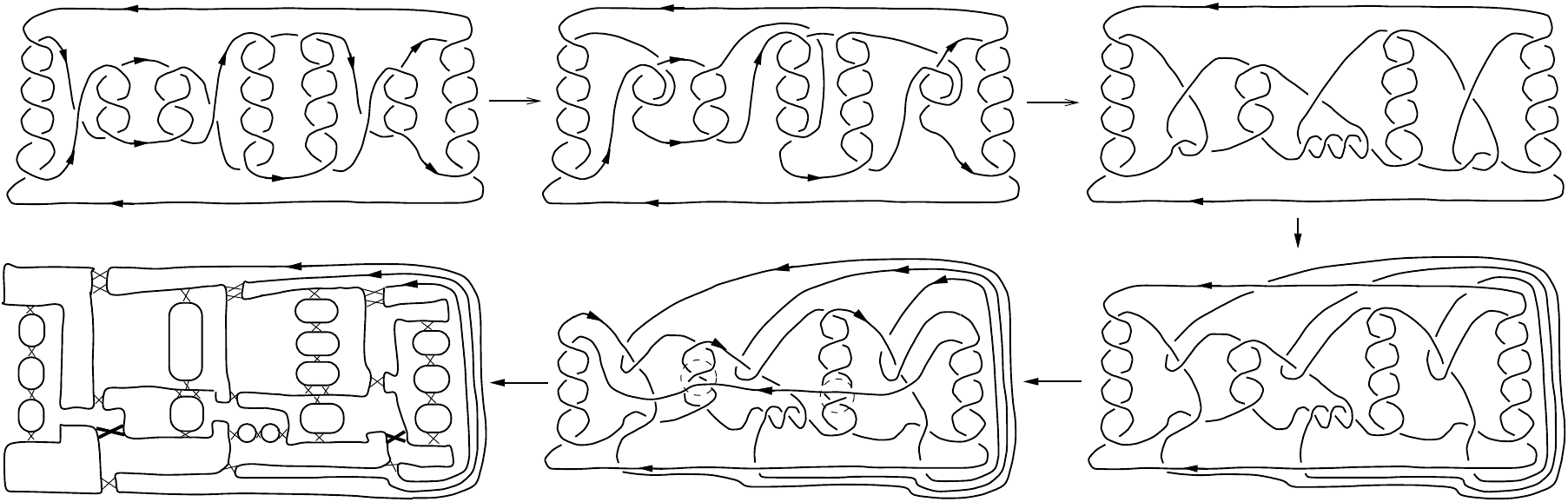}
\caption{The sequence of moves taking a diagram $D$ of $\L\in P_1(5,5,3,3;-5,-5,-3,-1,-1,-1)$ to a diagram $\tilde{D}$ such that $s(\tilde{D})=1-\kappa^++2\sum \alpha_j+2\sum \beta_i=19$. The two bold crossings indicate the horizontal strips corresponding to the two $\alpha_j$ strips containing 3 crossings. Notice that there are $-\kappa^++\sum \alpha_j+\sum \beta_i=7$ S/MP moves available hence $\L$ has a diagram with $12=1+\sum \alpha_j+\sum \beta_i$ Seifert circles. In this case, there are two compatible S moves and they must be performed on the under strands at the two crossings indicated by dashed circles.}
\label{beta_zero2}
\end{figure}

\medskip
Let us now consider the case when $\kappa^->\tau=\kappa^+-1$. Use cyclic permutation and flype moves, we can choose a standard diagram of $\L$ such that (i) The left most strip is a positive strip (say it is the $\alpha_1$ strip) and the right most strip is a negative strip with more than one crossing (say it is the $\beta_{\kappa^-}$ strip); (ii) for each of the remaining positive strips, place a negative single crossing strip to its left. See the first diagram in Figure \ref{beta_zero2}. The moves used in this case (as shown in Figure \ref{beta_zero2}) are similar to those used in for the case of $\tau=\kappa^-$, but with a slight modification. Here, for $2\le j\le \kappa^+$, each $\alpha_j$ strip is converted to a horizontal strip containing $2\alpha_j-1$ crossings, which contributes $2\alpha_j-2$ Seifert circles to the resulting diagram $\tilde{D}$. There are $\kappa^+-1=\tau$ concentric large Seifert circles in $\tilde{D}$, with the $\beta_{\kappa^-}$ strip, which now has $2\beta_{\kappa^-}$ crossings and contributes $2\beta_{\kappa^-}-1$ Seifert circles to $\tilde{D}$, contained in the inner most of these concentric large Seifert circles. For $1\le i\le \kappa^--1$ (in the case that $\kappa^->\kappa^+$), each $\beta_i$ strips still has the same number of crossings, is contained in one of these concentric large Seifert circles and has exactly one S move and $\beta_i-1$ MP moves. The under strands at the crossings where the S moves can be made are marked by dashed circles in Figure  \ref{beta_zero2}. Finally, the remaining $\alpha_1$ strip has $2\alpha_1$ crossings and contributes $2\alpha_1-1$ Seifert circles to $\tilde{D}$, and there are two additional Seifert circles not contained in the concentric large Seifert circles. Combining these, we see that $\tilde{D}$ has a total of $1-\kappa^++2\sum \alpha_j+2\sum \beta_i$ Seifert circles, with $-\kappa^++\sum\alpha_j+\sum\beta_i$ S/MP moves. It follows that $\L$ can be represented by a diagram with $1+\sum\alpha_j+\sum\beta_i$ Seifert circles.

\medskip
Proof for upper bound in (\ref{MT1e5}): This follows trivially from Remark \ref{repeated_MP_move}.

\medskip
We have now established that all the formulas given in Theorem \ref{MT1} are braid index upper bounds for the corresponding pretzel link diagrams.

\section{Braid index upper bounds for Type 2 pretzel links}\label{upperbound2}

\medskip
Proof for upper bound in  (\ref{MT2e1}): This follows trivially from Remark \ref{repeated_MP_move}.

\begin{figure}[htb!]
\includegraphics[scale=.6]{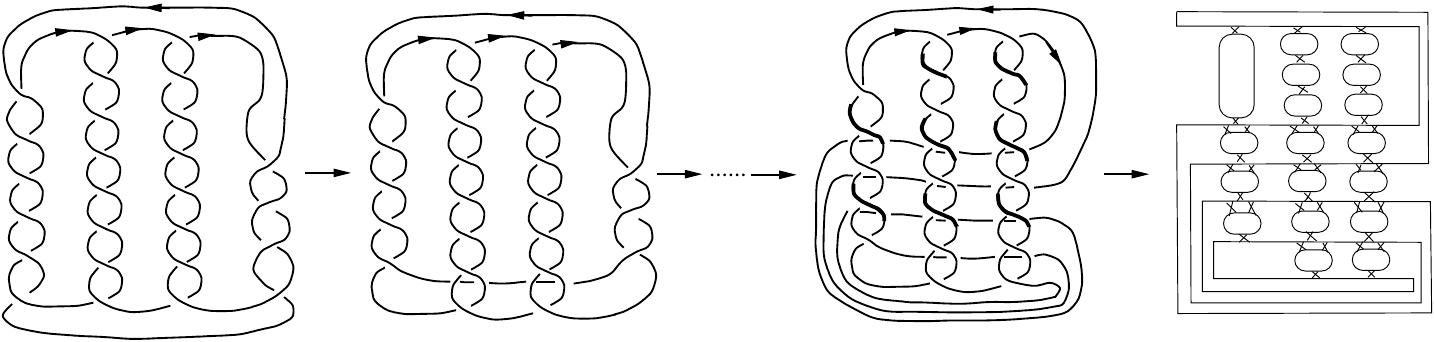}
\caption{The long moves taking the diagram $D$ of $\L\in P_2(6,8,8;-4)$ to $\tilde{D}$ with $s(\tilde{D})=19$ and 8 additional (compatible) S/MP moves.}
\label{beta1small2}
\end{figure}

\medskip
Proof for upper bound in  (\ref{MT2e2}): A sequence of long moves similar to the Type 1 (\ref{MT1e2}) and (\ref{MT1e3}) cases takes a standard diagram $D$ of $\L\in P_2(2\alpha_1,\ldots,2\alpha_{\kappa^+};-2\beta_1)$ to a new diagram $\tilde{D}$ as shown in Figure \ref{beta1small2} 
(after a reduction by $1+2\beta_1$ Seifert circles using $2\beta_1$ long moves). These moves eliminate the Seifert circles contributed by the $\beta_1$ strip and one Seifert circle contributed by the $\alpha_1$ strip. The two Seifert circles containing the top and bottom strands of $D$ are combined into one. Thus $s(\tilde{D})=-\kappa^++2\sum \alpha_j$. Furthermore, each $\alpha_j$ strip still contains $\alpha_j-1$ S/MP moves ($\beta_1$ S moves and $\alpha_j-\beta_1-1$ MP moves to be precise). One can verify that these S moves are compatible so the total number of realizable S/MP moves in $\tilde{D}$ is $-\kappa^++\sum \alpha_j$. It follows that $\L$ can be represented by a closed braid with $-\kappa^++2\sum \alpha_j-(-\kappa^++\sum \alpha_j)=\sum\alpha_j$ strings, as desired. This settles the case of (\ref{MT2e2}).

\begin{figure}[htb!]
\includegraphics[scale=.6]{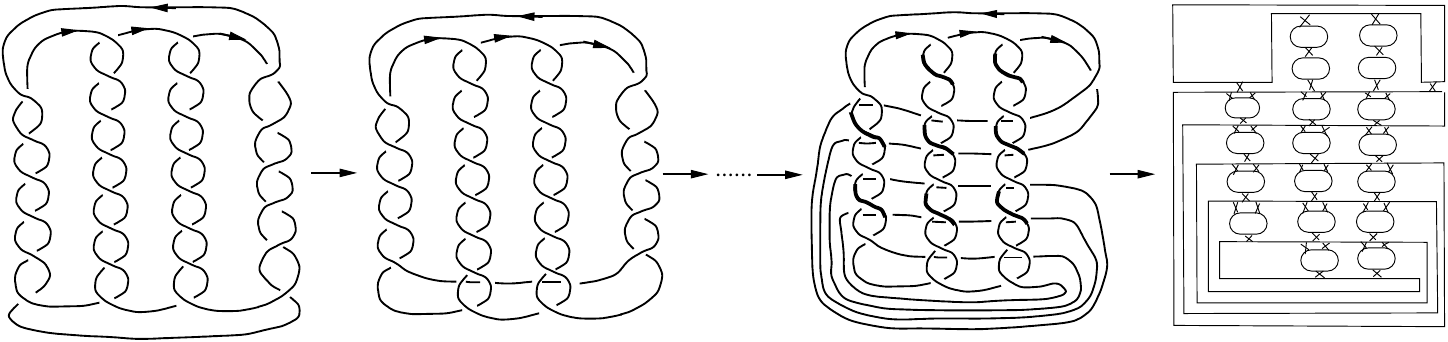}
\caption{The long moves taking the diagram $D$ of $\L\in P_2(6,8,8;-6)$ to $\tilde{D}$ with $s(\tilde{D})=20$ and 8 additional (compatible) S/MP moves.}
\label{beta1equal2}
\end{figure}

 \medskip
Proof for upper bound in  (\ref{MT2e3}): The same sequence of long moves used for the (\ref{MT2e2}) case takes a standard diagram $D$ of $\L\in P_2(2\alpha_1,\ldots,2\alpha_{\kappa^+};-2\beta_1)$ to a new diagram $\tilde{D}$ as shown in Figure \ref{beta1equal2} (after a reduction by $2\beta_1$ Seifert circles using $2\beta_1-1$ long moves). These moves eliminate the Seifert circles contributed by the $\beta_1$ strip and one Seifert circle contributed by the $\alpha_1$ strip. But in this case the two Seifert circles containing the top and bottom strands of $D$ remain as different Seifert circles (one of which is the long Seifert circle to be used for the S moves). Thus $s(\tilde{D})=-\kappa^++1+2\sum \alpha_j$. As before, each $\alpha_j$ strip still contains $\alpha_j-1$ S/MP moves ($\beta_1$ S moves and $\alpha_j-\beta_1-1$ MP moves). These S moves are compatible and it follows that $\L$ can be represented by a closed braid with $-\kappa^++1+2\sum \alpha_j-(-\kappa^++\sum \alpha_j)=1+\sum\alpha_j$ strings, as desired. In the case that $\beta_1>\min_{1\le j\le \kappa^+}\{\alpha_j\}$, the same discussion applies with the role of $\beta_1$ and $\alpha_{\kappa^+}$ switched (assuming that $\alpha_{\kappa^+}=\min_{1\le j\le \kappa^+}\{\alpha_j\}$), which leads to the desired result.
 This settles the case of (\ref{MT2e3}).
 
\begin{figure}[htb!]
\includegraphics[scale=.6]{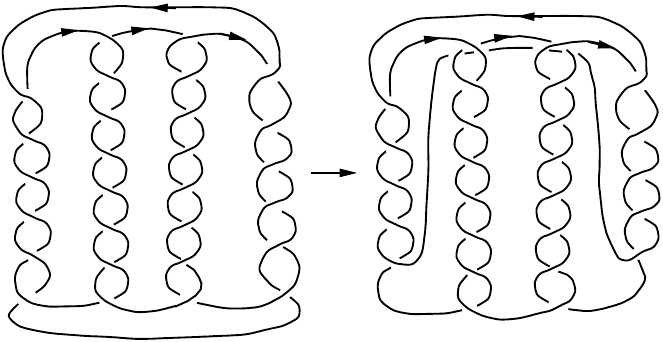}
\caption{The move taking the diagram $D$ of $\L\in P_2(6,8;-8,-6)$ to $\tilde{D}$ with $s(\tilde{D})=24=-\kappa^++2\sum \alpha_j-\kappa^-+2\sum \beta_i$ and $10=-\kappa^++\sum\alpha_j-\kappa^-+\sum\beta_i$ MP moves. Thus $\L$ can be represented by a diagram with $14=\sum \alpha_j+\sum \beta_i$ Seifert circles.}
\label{MT2e4fig}
\end{figure}

\medskip
Proof for upper bound in  (\ref{MT2e4}): Let  $\L\in P_{2}(2\alpha_1,\ldots, 2\alpha_{\kappa^+}; -2\beta_1,\ldots, -2\beta_{\kappa^-})$ and $D$ a standard diagram of $\L$, then the one move as shown in Figure \ref{MT2e4fig} reduces the total number of Seifert circles in $S(D)$ by two while keeping the number of MP moves unchanged. That is, there are still $-\kappa^++\sum\alpha_j-\kappa^-+\sum\beta_i$ MP moves hence
$$-\kappa^++2\sum\alpha_j-\kappa^-+2\sum\beta_i-(-\kappa^++\sum\alpha_j-\kappa^-+\sum\beta_i)=\sum \alpha_j+\sum \beta_i$$
is an upper bound of $\b(\L)$. Notice that the construction as shown in Figure \ref{MT2e4fig} is always possible since there must two adjacent strips with different signs. We have now established that all the formulas given in Theorem \ref{MT2} are braid index upper bounds for the corresponding pretzel link diagrams.
 
\section{Ending Remarks}

We end our paper with the following remark.

\medskip
\begin{remark}{\em 
In the case of Type 1 and Type 2 pretzel links, we always have $\b(\L)=\b_0(\L)$. We see that in some cases, the difference between the braid indices of a non-alternating pretzel link and its alternating counterpart can be as large as one wants. In the sequel of this paper, we shall complete our study on Type 3 pretzel links. We have found that it is not always true that $\b(\L)=\b_0(\L)$ for Type 3 pretzel links, even though it is still true for most of them. We are able to identify all such pretzel links for which it is possible that $\b(\L)>\b_0(\L)$. In these potential exceptional cases, we have further established that $\b_0(\L)\le \b(\L)\le 1+\b_0(\L)$.
}
\end{remark}

\end{document}